\documentclass[a4paper,11pt,fleqn]{article}
\usepackage[utf8]{inputenc}
\usepackage[T1]{fontenc}    
\usepackage[ruled,vlined]{algorithm2e}
\usepackage{amsmath,amsfonts,amssymb,amsthm}
\usepackage{bm,bbm}
\usepackage{mathtools}
\usepackage{mathrsfs} 
\usepackage{thm-restate,thmtools}
\usepackage{soul}
\usepackage{algorithmic}
\usepackage{nicefrac}
\usepackage{color}
\usepackage{xcolor,colortbl}
\usepackage{nicefrac}  
\usepackage{booktabs}
\usepackage{adjustbox}
\usepackage{hyperref}
\usepackage{url}
\usepackage{cleveref}
\usepackage[flushleft]{threeparttable}
\usepackage{multirow}
\usepackage[shortlabels]{enumitem}
\usepackage{pifont}
\usepackage{microtype}      
\usepackage{natbib}
\usepackage{doi}
\usepackage{authblk}
\usepackage{tikz}
\usepackage{tablefootnote}
\usepackage{needspace}
\usepackage{euscript}
\usepackage{charter}
\usepackage{array}
\usepackage[normalem]{ulem}
\usepackage{lipsum}
\usepackage{graphicx}
\graphicspath{ {./images/} }

\newtheorem{theorem}{Theorem}[section]
\newtheorem{proposition}[theorem]{Proposition}
\newtheorem{lemma}[theorem]{Lemma}

\newtheorem{definition}[theorem]{Definition}
\newtheorem{assumption}{Assumption}
\newtheorem{remark}[theorem]{Remark}

\definecolor{labelkey}{rgb}{0,0.08,0.45}
\definecolor{refkey}{rgb}{0,0.6,0.0}
\definecolor{Brown}{rgb}{0.45,0.0,0.05}
\definecolor{dgreen}{rgb}{0.00,0.49,0.00}
\definecolor{dblue}{rgb}{0,0.08,0.75}

\hypersetup{
colorlinks,
linkcolor=dgreen,
citecolor=dblue
}

\DeclareMathOperator*{\grad}{\mathit{\nabla \!f}}

\newcommand{\R}{\mathbb{R}}

\newcommand{\E}{\mathbb{E}}

\newcommand{\bbE}{\mathbb{E}}
\newcommand{\bbN}{\mathbb{N}}
\newcommand{\bbP}{\mathbb{P}}
\newcommand{\bbR}{\mathbb{R}}

\newcommand{\cP}{\mathcal{P}}

\newcommand{\cX}{\mathcal{X}}
\newcommand{\cZ}{\mathcal{Z}}

\newcommand{\dist}{\text{dist}}
\newcommand{\norm}[1]{\left\|#1\right\|}

\numberwithin{equation}{section}

\title{ { \sffamily Taming the Tails of Stochastic Gradient Descent\\
in the Interpolation Regime} } 

\author[1,3]{Leonardo Galli}
\author[2]{Daniela Angela Parletta}
\author[2]{Andrea Paudice}

\affil[1]{\footnotesize Mathematisches Institut, LMU Munich, Theresienstr. 39, 80333 Munich, Germany}
\affil[2]{\footnotesize Department of Computer Science, Aarhus University, Åbogade 34, 8200 Aarhus, Denmark}
\affil[3]{Munich Center for Machine Learning}

\date{}

\begin{document}
\maketitle

\begin{abstract}
In this paper, we study stochastic gradient descent in the smooth \emph{interpolation} regime, where the population objective shares a common minimizer with the sample losses. In contrast to prior work, which has largely focused on in-expectation analysis, we obtain high-probability convergence guarantees under the same assumptions. In the convex case, we prove a \emph{fast rate} of order $\log(\beta^{-1}) \cdot n^{-1}$ after $n$ iterations, with failure probability at most $\beta$. In the strongly convex case, we derive a family of moment-based bounds of order $\beta^{-p} \cdot \exp(-\mathcal{O}(p \cdot n))$ for $p \in (0,1]$, and complement them with a matching lower bound showing that the corresponding deterioration in the exponential rate is essentially unavoidable. In both settings, our theoretical results are unified across a broad class of step-size schedules. Building on this analysis, we extend our results to a low-noise regime and to nonconvex objectives satisfying the Polyak-\L{}ojasiewicz condition.
\end{abstract}

\vspace{1ex}
\allowdisplaybreaks

\section{Introduction}
Stochastic Gradient Descent (SGD) \cite{Robbins1951,Bottou2018} is among the most widely used algorithms in modern machine learning. Its performance is typically studied within the framework of stochastic optimization, where the goal is to minimize a population loss having access only to stochastic gradients, or equivalently to sample losses. The success of SGD stems from a combination of algorithmic simplicity, empirical effectiveness, and strong theoretical guarantees. For instance, on smooth objectives with stochastic gradients of bounded variance, SGD achieves, after $n$ iterations, rates of order $n^{-1/2}$ and $n^{-1}$ on convex and strongly convex objectives \cite{Nemirovski2009}, respectively. While optimal \cite{Nemirovski1983,Agarwal2012} in the classical stochastic setting, these rates are substantially slower than their deterministic counterparts, namely $n^{-1}$ and $\exp(- \mathcal{O}(n))$, which are achievable via Gradient Descent (GD) \cite{Nesterov2004}. 

A structural property that allows SGD to recover deterministic rates is \emph{interpolation}, which in its original formulation requires the population objective and the sample losses to share a common minimizer  \cite{Ma2018}. This regime is well motivated by modern over-parameterized learning models, where the number of parameters often far exceeds the number of training samples and perfect fitting of the training data is commonly observed. Interpolation has profound algorithmic consequences, and has been exploited to obtain SGD rates matching those of GD \cite{Ma2018,Vaswani2019b,Loizou2021,Galli2023}. Despite this progress, the existing theory for SGD under interpolation is largely confined to convergence \emph{in expectation}, which leaves unclear what can be guaranteed along a \emph{single}, potentially expensive, run of SGD. While in-expectation bounds can always be turned into high-probability guarantees via Markov's inequality, this typically yields a poor dependence on the failure probability $\beta \in (0,1)$. Obtaining substantially stronger high-probability rates is in general much more difficult, and often requires stronger assumptions (e.g., \cite{Nemirovski2009,Harvey2019,Liu2022,Eldowa2024,Liu2024}).

\paragraph{Contributions.}
In this paper, we show that, in the interpolation regime, one can do substantially better without strengthening the assumptions used in the corresponding in-expectation analysis. In particular, we provide the following contributions.

\begin{itemize}
\item For \emph{convex} smooth objectives, we establish, to the best of our knowledge, the first high-probability convergence rate that matches the corresponding in-expectation result up to a $\log(\beta^{-1})$ factor; see \Cref{thm:main_result}. Remarkably, the fast rate of order $n^{-1}$ is obtained without requiring assumptions beyond those used in the in-expectation analysis. Moreover, our analysis is unified across a broad class of step-size policies, well beyond the constant schedule. In particular, it covers Stochastic Line Search (SLS) \cite{Vaswani2019b} and Stochastic Polyak Step Size (SPS) \cite{Loizou2021}, adaptive methods that allow one to run SGD without additional knowledge of the underlying objective. This is made possible by exploiting the fact that, under interpolation, all these methods enjoy \emph{Fej\'er monotonicity} \cite{Fejer1922,Agmon1954}; see \Cref{prop:ss_condition_examples}. This property confines the iterates near the optimum and is what ultimately allows us to obtain the fast high-probability rate.

\item For \emph{strongly convex} objectives, we derive a family of high-probability bounds (\Cref{thm:markov_p}) of order $\beta^{-p}\cdot \exp(- \mathcal{O}(p \cdot \kappa^{-1} \cdot n))$ for $p \in (0,1]$, where $\kappa$ denotes the condition number of the population objective. The choice $p=1$ recovers the canonical Markov-type bound, which matches the exponential in-expectation rate in $n$ at the price of a poor dependence on $\beta$. Smaller values of $p$ improve the dependence on $\beta$, at the cost of a slower, yet still exponential, rate in $n$. We complement this result with a matching lower bound (\Cref{thm:sc_lower_bound_p}), showing that the deterioration of the exponential rate is essentially tight. The lower-bound construction reveals that this phenomenon can be traced to a ``lazy'' event on which all sampled gradients vanish.

\item We extend our analysis in two directions: \emph{beyond interpolation} and \emph{beyond convexity}. 
In \Cref{thm:hp_generalized_interpolation}, we consider a low-noise regime in which the gradient noise at the minimizer is subGaussian with modulus $\sigma_\star$. When $\sigma_\star = 0$, the fast $n^{-1}$ rate is recovered, whereas for $\sigma_\star> 0$ the rate degrades gracefully to order $n^{-1/2}$. This result relies crucially on the error decomposition developed for the interpolation regime and requires a careful high-probability control of the location of the iterates, which may grow as $\sqrt{n}$ in the absence of Fej\'er monotonicity. Finally, we move beyond convexity by considering functions satisfying the Polyak--\L ojasiewicz (PL) condition. In \Cref{thm:pl_qg_extension}, we prove an exponential high-probability rate of order $\beta^{-1/2}\cdot \exp(- \mathcal{O}(\kappa^{-2} \cdot n))$. In this setting, the lack of Fej\'er monotonicity is handled by exploiting the quadratic-growth property implied by the PL condition \cite{Karimi2016}.
\end{itemize}

\subsection{Related Works}
\noindent \textbf{SGD under Interpolation.}
The study of stochastic gradient methods in the interpolation regime was initiated in \cite{Ma2018}, where it was shown that SGD recovers, for strongly convex objectives, the same exponential rate as GD. Subsequent works established deterministic-like rates under interpolation for a variety of objectives and step-size rules. In particular, \cite{Vaswani2019a} considered the stronger Strong Growth Condition and proved deterministic rates for constant step-size SGD in the strongly convex, convex, nonconvex, and PL settings. Moving beyond fixed step-sizes, \cite{Vaswani2019b}, \cite{Loizou2021}, and \cite{Galli2023} analyzed several line-search and Polyak-type rules, including adaptive procedures that do not require prior knowledge of problem-dependent constants and often perform better in practice. Although the term interpolation originates from \cite{Ma2018}, the zero-noise-at-the-optimum specialization of the assumptions in \cite{Moulines2011} already encompasses this regime. Their rate, however, never matches $\mathcal{O}(n^{-1})$, even in the noiseless case where their parameter $\sigma$ is 0. More recently, last-iterate guarantees beyond exact interpolation have been investigated in \cite{Attia2025,Garrigos2025}; in particular, these works prove bounds of order $n^{-1} + \sigma_\star \cdot n^{-1/2}$, up to logarithmic factors. All these results are stated in expectation. On a different front, \cite{Sebbouh2021} established almost sure convergence rates for SGD and related methods in several convex and nonconvex settings. These guarantees, however, are asymptotic in nature and do not provide finite-sample bounds.

\noindent \textbf{High-Probability Bounds without Interpolation.}
Obtaining high-probability rates matching the corresponding in-expectation guarantees, while improving on the confidence dependence given by Markov's inequality, is generally difficult and typically requires stronger noise assumptions than the finite-variance condition used in in-expectation analysis. For instance, \cite{Harvey2019,Jain2021,Liu2023,Liu2024,Madden2024} derive high-probability convergence guarantees for SGD and related methods under subGaussian noise assumptions, yielding \emph{subGaussian} rates of order $\sqrt{\log(\beta^{-1}) \cdot n^{-1}}$, which are optimal in the absence of additional structure. To recover deterministic rates, \cite{Jin2024} imposes even stronger noise assumptions, including boundedness; in the convex setting, for instance, this effectively requires a bounded domain. In \cite{Bach2013}, a high-probability fast rate is derived for the specific case of finite-dimensional least-squares regression under almost-sure boundedness of the covariates, higher-order moment control of the residual noise, and a uniform bound on the directional kurtosis of the covariates. A complementary line of work considers robust variants based on clipping or truncation, which can recover sharp rates under assumptions similar to those used for in-expectation guarantees; representative examples include \cite{Gorbunov2020,Nguyen2023,Parletta2024}. Along a different direction, \cite{Eldowa2024,Lou2022} study more general noise models and prove a two-regime behavior for SMD, and therefore also for SGD, in which the contribution of the heavier-tailed component decays at rate $n^{-1}$; as a consequence, the asymptotic behavior is eventually governed by the subGaussian term.

\section{Problem Setting}\label{sec:problem_setting}
We consider the following stochastic optimization problem
\[
\underset{x \in \mathcal{X}}{\text{minimize}} \quad F(x) = \mathbb{E}_{Z \sim \cP}[f(x, Z)]
\]
where $Z$ is a $\mathcal{Z}$-valued random variable with distribution $\cP$. 

\noindent We make the following assumptions.
\begin{assumption}[Feasibility Domain]
\label{ass:feasibility}
$\cX$ denotes a real separable Hilbert space, and $\cX_{\star} \coloneqq \arg \min_{x \in \mathcal{X}} F(x) \neq \emptyset$. We also define  $F_{\star}=\min_{x \in \cX} F(x) \in\mathbb{R}$.
\end{assumption}
\begin{assumption}[Convexity and Smoothness]
\label{ass:convexity}
For $\cP$-almost every $z$, $f(\cdot, z)\colon \mathcal{X} \rightarrow \mathbb{R}$ is convex and $L$-smooth\footnote{This means that for $\cP$-almost every $z$, $\nabla f(\cdot,z)$ is $L$-Lipschitz.} (this implies that $F$ is convex and $L$-smooth).
\end{assumption}
\begin{assumption}[Interpolation Regime]
\label{ass:interpolation}
There exists $x_{\star} \in \cX_{\star}$ s.t. $x_{\star} \in \arg\min_{x \in \cX} f(x, z)$ $\cP$-almost every $z$.
\end{assumption}
\paragraph{Notation.} Throughout $k$ and $n$ will always denote natural numbers. For $m \in \bbN$, let us define $[m] \coloneqq \{1,\dots,m\}$. As we consider iterative algorithms, we will denote by $Z_k \sim \cP$ the sample used by the algorithm at iteration $k$. For each $k \in \mathbb{N}$, we denote by $\mathbb{E}_k[ \cdot ] \coloneqq \mathbb{E}[ \cdot | \sigma(Z_1,\dots,Z_{k-1})]$. Note that
\ifdefined\itistwocolumns
\begin{equation}\label{eq:cond_exp}
\begin{split}
\E_k[f(x_k, Z_k)] &= F(x_k), \\
\E_k[\nabla f(x_k, Z_k)] &= \nabla F(x_k). 
\end{split}
\end{equation}
\else
\begin{equation}\label{eq:cond_exp}
\E_k[f(x_k, Z_k)] = F(x_k) \quad \text{and} \quad \E_k[\nabla f(x_k, Z_k)] = \nabla F(x_k).
\end{equation}
\fi
\noindent For $x \in \mathcal{X}$ and $z \in \cZ$, let us define $\Delta(x) \coloneqq F(x)-F_{\star}$, and (with a slight abuse of notation) $\Delta(x, z) \coloneqq f(x, z) - f(x_{\star}, z)$. 
\paragraph{Stochastic Gradient Descent.} Let $x_1 \in \mathcal{X}$ and $(\eta_k)_{k \geq 1}$ be a (possibly random) sequence of non-negative reals. SGD is defined by the following recurrence
\begin{align}
\label{eq:sgd_ls}
\forall k \geq 1, \quad x_{k+1} = x_k - \eta_k \cdot \grad(x_k,Z_k).
\end{align}
For a given $n \geq 1$, alongside $x_{n+1}$ we will also be interested in the average iterate
\begin{align*}
\bar x_n = \frac 1 n \sum_{k=1}^n x_k.
\end{align*}
We will assume that \emph{step size} sequence $\eta_k$ satisfies the following assumption. 
\begin{assumption}[Step Size Condition]
\label{ass:ss_condition}
Let $x_{\star} \in \cX_{\star}$. There exists a (possibly randomized) sequence $(\phi_k)_{k\geq 1}$ and a constant $\bar{\phi}>0$ s.t., for every $k \geq 1$, it holds that: $\phi_k \geq \bar{\phi}$ and \ifdefined\itistwocolumns (almost surely)
\begin{align}
\label{eq:ss_condition}
\|x_{k+1}-x_{\star}\|^2
\leq
\|x_k-x_{\star}\|^2-2 \cdot \phi_k\cdot \Delta(x_k,Z_k).
\end{align}
\else
\begin{align}
\label{eq:ss_condition}
\|x_{k+1}-x_{\star}\|^2
\leq
\|x_k-x_{\star}\|^2-2 \cdot \phi_k\cdot \Delta(x_k,Z_k) \quad \text{(almost surely)}.
\end{align}
\fi
\end{assumption}
\noindent The next proposition shows that, under interpolation, Assumption~\ref{ass:ss_condition} is satisfied by several step size policies of interest. The proof is in Appendix \ref{subsec:appendix_convex_smooth_optimization_prop_ss_condition_examples}.
\begin{proposition}\label{prop:ss_condition_examples}
Suppose Assumptions~\ref{ass:feasibility}-\ref{ass:interpolation} hold. Then Assumption~\ref{ass:ss_condition} is satisfied for the following methods:
\begin{enumerate}[label={\rm(\roman*)}]
\item\label{prop:ss_condition_eta} \textbf{Constant Step Size.}
If $\eta_k=\eta\in(0,2/L)$ for every $k\ge 1$, then Assumption~\ref{ass:ss_condition} holds with
\ifdefined\itistwocolumns
\begin{align}\label{eq:ss_condition_eta}
\phi_k =
\begin{cases}
\eta, & \text{if } \eta\in(0, 1/L],\\
\eta \cdot (2 - L \cdot \eta), & \text{if } \eta\in(1/L, 2/L),
\end{cases}    
\end{align}
and $\bar{\phi} = \phi_k$.
\else
\begin{align}\label{eq:ss_condition_eta}
\phi_k = \bar{\phi}  =
\begin{cases}
\eta, & \text{if } \eta\in(0, 1/L],\\
\eta \cdot (2 - L \cdot \eta), & \text{if } \eta\in(1/L, 2/L).
\end{cases}    
\end{align}
\fi

\item\label{prop:ss_condition_sls} \textbf{Stochastic Line Search (SLS).}
Suppose the step size $\eta_k$ is selected by the stochastic line-search rule in \Cref{alg:backtrack} (in Appendix \ref{sec:appendix_convex_smooth_optimization}), with parameters $c \in (1/2, 1)$, $\delta \in (0,1)$, and $\eta > 0$. Then Assumption~\ref{ass:ss_condition} holds with
\ifdefined\itistwocolumns
\begin{equation}\label{eq:ss_condition_sls}
    \begin{split}
\phi_k&=\frac{(2 \cdot c - 1)}{2 \cdot c} \cdot \eta_k,\\
\bar\phi&=\frac{(2 \cdot c - 1)}{2 \cdot c} \cdot \min\bigg\{\frac{2 \cdot (1-c) \cdot \delta}{L}, \eta\bigg\}.
    \end{split}
\end{equation}
\else
\begin{align}
\phi_k=\frac{(2 \cdot c - 1)}{2 \cdot c} \cdot \eta_k,
\quad
\bar\phi=\frac{(2 \cdot c - 1)}{2 \cdot c} \cdot \min\bigg\{\frac{2 \cdot (1-c) \cdot \delta}{L}, \eta\bigg\}.
\end{align}
\fi
\item\label{prop:ss_condition_polyak} \textbf{Stochastic Polyak Step size (SPS).}
Suppose
\begin{align*}
\eta_k=
\begin{cases}
\dfrac{\Delta(x_k, Z_k)}{\|\nabla f(x_k,Z_k)\|^2},
& \text{if } \|\nabla f(x_k,Z_k) \| \neq 0,\\[3mm]
0,
& \text{otherwise}.
\end{cases}
\end{align*}
Then Assumption~\ref{ass:ss_condition} holds with 
\begin{align}\label{eq:ss_condition_polyak}
\phi_k = \bar\phi = \frac{1}{4 \cdot L}.
\end{align}
\end{enumerate}
\end{proposition}
A few comments are in order. Assumption~\ref{ass:ss_condition} is introduced to isolate the one-step descent property that drives our analysis. In this way, it provides a unified treatment of different step-size policies and makes clear that our results immediately extend to any other schedule for which the condition can be verified. Among the examples covered by \Cref{prop:ss_condition_examples}, the constant step-size schedule is the canonical choice \cite{Nesterov2004}, and its best theoretical guarantees are typically obtained for $\eta_k=1/L$. In practice, however, the Lipschitz constant $L$ may be unknown or overly conservative, making such a choice infeasible or suboptimal. This motivates adaptive procedures such as those studied in \cite{Vaswani2019b,Loizou2021}, which do not require prior knowledge of $L$ and may also allow for larger step sizes \cite{Galli2023}. Finally, note that the SPS requires knowledge of $f(x_\star,Z_k)$ for every $k\geq 1$; this is often a mild prerequisite in applications where the minimum value of each sample loss is typically $0$.
\section{Convex Smooth Optimization}\label{sec:convex_smooth_optimization}
The following is the main result of this section.
\begin{theorem}
\label{thm:main_result}
\ifdefined\itistwocolumns
Suppose Assumptions~\ref{ass:feasibility}-\ref{ass:ss_condition} hold w.r.t. $x_{\star} \in \cX_{\star}$,  and let $x_1 \in \mathcal{X}$. Then, the iterates generated by \eqref{eq:sgd_ls} satisfy the following. For every $n \geq 1$ and  $\beta \in (0,1)$, with probability at least $1-\beta$, it holds that
\begin{equation}
\begin{split}
\label{eq:main_result}
\Delta(\bar{x}_n) \leq  \bigg[ \frac{1}{\bar{\phi} \cdot L} + \frac{4}{3} \cdot \log \Big( \frac{1}{\beta} \Big)\bigg]\cdot \frac{L \cdot \|x_1-x_{\star}\|^2}{n}.
\end{split}
\end{equation}
\else
Suppose Assumptions~\ref{ass:feasibility}-\ref{ass:ss_condition} hold w.r.t. $x_{\star} \in \cX_{\star}$,  and let $x_1 \in \mathcal{X}$. Then, the iterates generated by \eqref{eq:sgd_ls} satisfy the following. For every $n \geq 1$ and  $\beta \in (0,1)$, it holds that
\begin{align}
\label{eq:main_result}
\mathbb{P}\bigg(\Delta(\bar{x}_n) &\leq  \bigg[ \frac{1}{\bar{\phi} \cdot L} + \frac{4}{3} \cdot \log \Big( \frac{1}{\beta} \Big)\bigg] \cdot \frac{L \cdot \|x_1-x_{\star}\|^2}{n} \bigg) \geq 1 - \beta.
\end{align}
\fi
\end{theorem}
State-of-the-art in-expectation guarantees for SGD in this setting are of order $L \cdot \|x_1-x_{\star}\|^2 \cdot n^{-1}$ \cite{Vaswani2019b,Loizou2021}. 
Compared with these results, \Cref{thm:main_result} only incurs a $\log(\beta^{-1})$ overhead. Remarkably, unlike existing high-probability analyses in related settings, this fast rate is obtained without requiring assumptions beyond those used for the corresponding in-expectation guarantee. Moreover, although \Cref{thm:main_result} is formulated through the abstract \Cref{ass:ss_condition}, interpolation remains the key structural property in the present framework, as it is precisely what ensures that \Cref{ass:ss_condition} holds for several concrete step-size policies of interest, including constant step size, SLS, and SPS.

\begin{remark}
\label{rmk:main_result}
In statistical learning, $\Delta(x)$ is often referred to as the \emph{excess risk}, and rates of order $n^{-1}$ are called \emph{fast rates}. Obtaining fast rates with high probability is notoriously difficult. Classical examples include, among others, Empirical Risk Minimization (ERM) in realizable binary classification \cite{Blumer1989} and SGD on strongly convex risks (e.g. \cite{Harvey2019}). \Cref{thm:main_result} complements existing results by showing that high-probability fast rates can also be obtained for SGD in the smooth interpolation regime. The result closest in spirit to ours is \cite{Klochkov2021}, where the authors show that gradient descent on the empirical risk achieves a rate of order $\log(\beta^{-1}) \cdot \log(n) \cdot n^{-1}$ under the following \emph{Bernstein condition} \cite{Bernstein1924}: there exists $B>0$ such that
\begin{align}
\label{eq:bernstein_condition_original}
\forall x\in \mathcal{X}, \qquad \bbE[\Delta(x, Z)^2] \leq B \cdot \Delta(x).
\end{align}
We note that our assumptions do yield a related bound along the SGD trajectory, namely
\[
\bbE_k\!\left[(\Delta(x_k,Z_k)-\Delta(x_k))^2\right] \leq \frac{L \cdot \|x_1-x_{\star}\|^2}{2} \cdot \Delta(x_k),
\]
see \Cref{lemma:mds_prop}~\cref{lemma:mds_conditional_variance}. However, this should not be confused with eq.~\eqref{eq:bernstein_condition_original}, and in particular it does not imply it. Indeed, eq.~\eqref{eq:bernstein_condition_original} is a global condition that must hold for every $x\in\cX$, whereas our bound is local in the sense that it is only required along the random iterates generated by the algorithm.
\end{remark}
\subsection{Proof of Theorem \ref{thm:main_result}}\label{sub_section:proof_main_result}
Here we prove \Cref{thm:main_result}, highlighting the main ideas. The proofs of Lemmas~\ref{lemma:mds_prop} and \ref{lemma:bernstein_condition} are deferred to Appendix~\ref{sec:appendix_convex_smooth_optimization}.

At a high level, the main difficulty is that high-probability bounds in stochastic optimization typically incur deviations of order $\sqrt{\log(\beta^{-1})/n}$, due to the martingale fluctuations generated by the sampling noise. A direct concentration argument would therefore not yield the fast rate in \Cref{thm:main_result}. The key point is that, under \Cref{ass:ss_condition}, the iterates enjoy Fej\'er monotonicity with respect to $x_\star$, and are therefore confined to a deterministic ball around the optimum. This makes the relevant martingale self-normalized, in the sense that its conditional variance can be bounded in terms of the optimization errors $(\Delta(x_k))_{k\geq 1}$. As a result, the martingale term can be controlled by two contributions: one proportional to $\sum_{k=1}^n \Delta(x_k)$, which can be absorbed into the left-hand side, and one proportional to $\log(1/\beta)$ times the squared initial distance, which yields the final $n^{-1}$ rate after averaging.

\paragraph{Bounding Empirical Errors.} For $k \geq 1$, thanks to \Cref{ass:ss_condition}, it holds (almost surely) that
\begin{align}
\label{eq:one_step_inequality}
\|x_{k+1}-x_{\star}\|^2 &\leq \|x_{k} -x_{\star}\|^2 - 2 \cdot \phi_k \cdot \Delta(x_k, Z_k),
\end{align}
Rearranging and summing all the inequalities from $1$ to $n$, yields
\ifdefined\itistwocolumns
\begin{align*}
&\sum_{k=1}^n \phi_k \cdot (f(x_k, Z_k) - f(x_{\star}, Z_k))\leq\\
&\frac{1} {2} \cdot \Bigg(\|x_{1} -x_{\star}\|^2 - \|x_{n+1}-x_{\star}\|^2\Bigg) \leq \nonumber\frac{\|x_{1} -x_{\star}\|^2}{2} \cdot 
\end{align*}
\else
\begin{align*}
\sum_{k=1}^n \phi_k \cdot (f(x_k, Z_k) - f(x_{\star}, Z_k))  &\leq \frac{1} {2} \cdot \Bigg(\|x_{1} -x_{\star}\|^2 - \|x_{n+1}-x_{\star}\|^2\Bigg) \nonumber\\
&\leq \frac{\|x_{1} -x_{\star}\|^2}{2} \cdot 
\end{align*}
\fi
Using \Cref{ass:ss_condition} on the l.h.s. and dividing both sides of the above display by $n \cdot \bar \phi$ yields 
\begin{equation*}
\begin{split}
\frac{1}{n} \cdot \sum_{k=1}^n \Delta(x_k, Z_k)  &\leq \frac{\|x_{1} -x_{\star}\|^2}{2 \cdot \bar{\phi} \cdot n}. 
\end{split}
\end{equation*}
Then 
\ifdefined\itistwocolumns
\begin{equation}
\label{eq:as_bound}
\begin{split}
\frac{1}{n} \cdot \sum_{k=1}^n \Delta(x_k) &= \frac{1}{n} \cdot \sum_{k=1}^n \Delta(x_k, Z_k) \\
&\quad + \frac{1}{n} \cdot \sum_{k=1}^n \underbrace{(\Delta(x_k) - \Delta(x_k, Z_k))}_{\zeta_k}\\
&\leq \frac{\|x_{1} -x_{\star}\|^2}{2 \cdot \bar{\phi} \cdot n} + \frac{1}{n} \cdot \underbrace{\sum_{k=1}^n \zeta_k.}_{G_n}
\end{split}
\end{equation}
\else
\begin{align}
\label{eq:as_bound}
\frac{1}{n} \cdot \sum_{k=1}^n \Delta(x_k) &= \frac{1}{n} \cdot \sum_{k=1}^n \Delta(x_k, Z_k) + \frac{1}{n} \cdot \sum_{k=1}^n \underbrace{(\Delta(x_k) - \Delta(x_k, Z_k))}_{\zeta_k} \nonumber\\
&\leq \frac{\|x_{1} -x_{\star}\|^2}{2 \cdot \bar{\phi} \cdot n} + \frac{1}{n} \cdot \underbrace{\sum_{k=1}^n \zeta_k.}_{G_n}
\end{align}
\fi
Note also by eq.~\eqref{eq:one_step_inequality},  the sequence $(x_{k})_{k \in \mathbb{N}}$ enjoys \emph{Fej\'er monotonicity} w.r.t. $x_{\star}$ \cite{Fejer1922,Agmon1954}, i.e.
\begin{align}
\label{eq:fejer_monotonicity}
\|x_{k+1}-x_{\star}\| \leq \|x_k - x_{\star} \| \quad \text{(almost surely)}
\end{align}
which implies $\|x_k-x_{\star}\| \leq \|x_1-x_{\star}\| $ for every $k \geq 1$.

\paragraph{Bounding $G_n$.} We start with the following that establishes the relevant properties of the sequence $(\zeta_k)_{k \in [n]}$.
\begin{lemma}
\label{lemma:mds_prop}
Let $(\zeta_k)_{k \geq 1}$ be defined as in eq.~\eqref{eq:as_bound}, and suppose Assumptions \ref{ass:feasibility}-\ref{ass:ss_condition} hold. Then the following statements hold:
\begin{enumerate}[label={\rm(\roman*)}]
\item\label{lemma:mds_boundedness} For every $k \geq 1$, $\zeta_k$ satisfies \ifdefined\itistwocolumns (almost surely)
\begin{align}
\label{eq:mds_boundedness}
|\zeta_k| &\leq \frac{L \cdot \|x_1 - x_{\star}\|^2}{2}.
\end{align}
\else
\begin{align}
\label{eq:mds_boundedness}
|\zeta_k| &\leq \frac{L \cdot \|x_1 - x_{\star}\|^2}{2}  \quad \text{(almost surely)}
\end{align}
\fi
\item\label{lemma:mds} $(\zeta_k)_{k \in \bbN}$ is a martingale difference sequence.
\item\label{lemma:mds_conditional_variance} For every $k \geq 1$, $\zeta_k$ satisfies \ifdefined\itistwocolumns (almost surely)
\begin{align}
\label{eq:mds_conditional_variance}
\mathbb{E}_k[\zeta_k^2] &\leq \frac{L \cdot \|x_1 - x_{\star}\|^2}{2} \cdot \Delta(x_k).
\end{align}
\else
\begin{align}
\label{eq:mds_conditional_variance}
\mathbb{E}_k[\zeta_k^2] &\leq \frac{L \cdot \|x_1 - x_{\star}\|^2}{2} \cdot \Delta(x_k) \quad \text{(almost surely).}
\end{align}
\fi
\end{enumerate}
\end{lemma}
\noindent Next, we show that \cref{lemma:mds_boundedness}, \cref{lemma:mds}, and \cref{lemma:mds_conditional_variance} in \Cref{lemma:mds_prop} imply a (conditional) Bernstein's condition for $(\zeta_k)_{k \in \bbN}$.
\begin{lemma}[Bernstein's Condition]
\label{lemma:bernstein_condition}
Let $(\zeta_k)_{k \in [n]}$ be defined as in eq.~\eqref{eq:as_bound}, and suppose Assumptions \ref{ass:feasibility}-\ref{ass:ss_condition} hold. Then, for every $k \in [n]$ and for every $\lambda \in [0, 3/\nu)$, it holds (almost surely) that
\begin{align}
\label{eq:bernstein_condition}
\bbE_k[\exp(\lambda \cdot \zeta_k)] \leq \exp\bigg(\frac{\lambda^2 \cdot \nu }{2\cdot (1-\lambda \cdot \nu/3)} \cdot \Delta(x_k) \bigg),
\end{align}
where $\nu = \frac{L \cdot \|x_1-x_{\star}\|^2}{2}$.
\end{lemma}
\noindent In what follows, we denote by $\nu$ the constant defined in \Cref{lemma:bernstein_condition}. For $\lambda \in [0, 3/\nu)$, let $(\Psi_k(\lambda))_{k \geq 0}$ be defined as $\Psi_0(\lambda) \coloneqq 1$ and \ifdefined\itistwocolumns $\forall k \geq 1,$
\[
\Psi_k(\lambda) \coloneqq \exp\bigg(\lambda \cdot G_k - \frac{\lambda^2 \cdot \nu }{2 \cdot (1-\lambda \cdot \nu/3)} \sum_{i=1}^k \Delta(x_i)\bigg).
\]
\else
\[
\forall k \geq 1, \quad \Psi_k(\lambda) \coloneqq \exp\bigg(\lambda \cdot G_k - \frac{\lambda^2 \cdot \nu }{2 \cdot (1-\lambda \cdot \nu/3)} \sum_{i=1}^k \Delta(x_i)\bigg).
\]
\fi
Note that $\Psi_k(\lambda)$ is a supermartingale for every $\lambda \in [0,3/\nu)$. Indeed, for each $k \geq 1$ it holds that
\ifdefined\itistwocolumns
\begin{align*}
&\mathbb{E}_k[\Psi_k(\lambda)] =\\
&\mathbb{E}_k\Bigg[\exp\bigg(\lambda \cdot G_k - \frac{\lambda^2 \cdot \nu }{2 \cdot (1-\lambda \cdot \nu/3)} \cdot\sum_{i=1}^k \Delta(x_i)\bigg)\Bigg] =\\
&\mathbb{E}_k\Bigg[\exp\bigg(\lambda \cdot \sum_{i=1}^k \zeta_i\bigg)\Bigg] \cdot \\
&\quad \exp\bigg(- \frac{\lambda^2 \cdot \nu }{2 \cdot (1-\lambda \cdot \nu/3)} \cdot \sum_{i=1}^k \Delta(x_i)\bigg) \overset{(a)}{=}\\
&\mathbb{E}_k\Bigg[\exp\bigg(\lambda \cdot \zeta_k\bigg)\Bigg] \cdot \\
&\quad \exp\bigg( \lambda \cdot \sum_{i=1}^{k-1} \zeta_i - \frac{\lambda^2 \cdot \nu }{2 \cdot (1-\lambda \cdot \nu/3)} \cdot \sum_{i=1}^k \Delta(x_i)\bigg) \overset{(b)}{\leq}\\
&\exp\bigg(\frac{\lambda^2 \cdot \nu}{2 \cdot (1-\lambda \cdot \nu/3)} \cdot \Delta(x_k)\bigg)\cdot \\
&\quad\exp\bigg( \lambda \cdot \sum_{i=1}^{k-1} \zeta_i - \frac{\lambda^2 \cdot \nu }{2 \cdot (1-\lambda \cdot \nu/3)} \cdot \sum_{i=1}^k \Delta(x_i)\bigg) =\\
&\exp\bigg( \lambda \cdot \sum_{i=1}^{k-1} \zeta_i - \frac{\lambda^2 \cdot \nu }{2 \cdot (1-\lambda \cdot \nu/3)} \cdot \sum_{i=1}^{k-1} \Delta(x_i)\bigg) =\\
&\Psi_{k-1}(\lambda),
\end{align*}
\else
\begin{align*}
\mathbb{E}_k[\Psi_k(\lambda)] &= \mathbb{E}_k\Bigg[\exp\bigg(\lambda \cdot G_k - \frac{\lambda^2 \cdot \nu }{2 \cdot (1-\lambda \cdot \nu/3)} \cdot \sum_{i=1}^k \Delta(x_i)\bigg)\Bigg] \\
&= \mathbb{E}_k\Bigg[\exp\bigg(\lambda \cdot \sum_{i=1}^k \zeta_i\bigg)\Bigg]\exp\bigg(- \frac{\lambda^2 \cdot \nu }{2 \cdot (1-\lambda \cdot \nu/3)} \cdot \sum_{i=1}^k \Delta(x_i)\bigg) \\
&\overset{(a)}{=} \mathbb{E}_k\Bigg[\exp\bigg(\lambda \cdot \zeta_k\bigg)\Bigg]\exp\bigg( \lambda \cdot \sum_{i=1}^{k-1} \zeta_i - \frac{\lambda^2 \cdot \nu }{2 \cdot (1-\lambda \cdot \nu/3)} \cdot \sum_{i=1}^k \Delta(x_i)\bigg) \\
&\overset{(b)}{\leq} \exp\bigg(\frac{\lambda^2 \cdot \nu}{2 \cdot (1-\lambda \cdot \nu/3)} \cdot \Delta(x_k)\bigg)\exp\bigg( \lambda \cdot \sum_{i=1}^{k-1} \zeta_i - \frac{\lambda^2 \cdot \nu }{2 \cdot (1-\lambda \cdot \nu/3)} \cdot \sum_{i=1}^k \Delta(x_i)\bigg) \\
&= \exp\bigg( \lambda \cdot \sum_{i=1}^{k-1} \zeta_i - \frac{\lambda^2 \cdot \nu }{2 \cdot (1-\lambda \cdot \nu/3)} \cdot \sum_{i=1}^{k-1} \Delta(x_i)\bigg) \\
&= \Psi_{k-1}(\lambda),
\end{align*}
\fi
where (a) follows since $\zeta_i$, for each $i \leq k-1$, is $\sigma(Z_1,\dots,Z_{k-1})$-measurable, and (b) follows from \Cref{lemma:bernstein_condition}. 
Thus, by the tower rule of expectation, $\mathbb{E}[\Psi_n(\lambda)] \leq 1$. Note that, for every $t > 0$, it holds that
\ifdefined\itistwocolumns
\begin{equation}\label{eq:Gn_bound1}
    \begin{split}
&\mathbb{P}\bigg(G_n - \frac{\lambda \cdot \nu }{2 \cdot (1-\lambda \cdot \nu/3)} \cdot \sum_{k=1}^n \Delta(x_k) > t \bigg) = \\
&\mathbb{P}\big(\Psi_n(\lambda) > \exp(\lambda \cdot t) \big)\overset{(a)}{\leq} \\
& \exp(-\lambda \cdot t) \cdot \mathbb{E}[\Psi_n(\lambda)] \overset{(b)}{\leq} \exp(-\lambda \cdot t),
    \end{split}
\end{equation}
\else
\begin{align}
\label{eq:Gn_bound1}
\mathbb{P}\bigg(G_n - \frac{\lambda \cdot \nu }{2 \cdot (1-\lambda \cdot \nu/3)} \cdot \sum_{k=1}^n \Delta(x_k) > t \bigg) &= \mathbb{P}\big(\Psi_n(\lambda) > \exp(\lambda \cdot t) \big) \nonumber\\
&\overset{(a)}{\leq} \exp(-\lambda \cdot t) \cdot \mathbb{E}[\Psi_n(\lambda)] \nonumber\\
&\overset{(b)}{\leq} \exp(-\lambda \cdot t),
\end{align}
\fi
where (a) follows from Markov inequality, and (b) follows from the tower rule of expectation and $(\Psi_k(\lambda))_{k \geq 0}$ being a supermartingale for every $\lambda \in [0,3/\nu)$. Now, for $\beta \in (0,1)$, let $t = \log(1/\beta)/\lambda$. Then, eq.~\eqref{eq:Gn_bound1} implies that with probability at least $1-\beta$
\begin{align}
\label{eq:Gn_bound2}
G_n  \leq \frac{\lambda \cdot \nu }{2 \cdot (1-\lambda \cdot \nu/3)} \cdot\sum_{k=1}^n \Delta(x_k) + \frac{\log(1/\beta)}{\lambda}.
\end{align}
Choosing $\lambda = 3/(4\nu) \in [0, 3/\nu)$ yields
\begin{align}
\label{eq:lambda_choice}
\frac{\lambda \nu }{2(1-\lambda \nu/3)} = \frac{1}{2}.
\end{align}
Therefore, combining eq.~\eqref{eq:lambda_choice} with eq.~\eqref{eq:Gn_bound2} yields
\begin{align}
\label{eq:Gn_final_bound}
G_n  \leq \frac{1}{2} \cdot \sum_{k=1}^n \Delta(x_k) + \frac{4}{3} \cdot \nu \cdot \log\bigg(\frac 1 \beta \bigg).
\end{align}
\paragraph{Final Bound.} The thesis follows by plugging eq.~\eqref{eq:Gn_final_bound} in eq.~\eqref{eq:as_bound}, rearranging, using Jensen's inequality, and the definition of $\nu$.
\section{Strongly Convex Smooth Optimization}\label{sec:sc}
In this section, we turn to the strongly convex setting and consider the following. All proofs can be found in Appendix \ref{sec:appendix_sc}.
\begin{assumption}[Strong Convexity and Smoothness]
\label{ass:strong_convexity}
We assume that $F\colon\mathcal{X} \rightarrow \bbR$ is strongly convex with modulus $\mu > 0$, and that for $\cP$-almost every $z \in \cZ$, $f(\cdot, z)\colon \mathcal{X} \rightarrow \mathbb{R}$ is convex and $L$-smooth (this implies that $F$ is $L$-smooth). We define the \emph{stochastic condition number} as $\kappa \coloneqq L/\mu$.
\end{assumption}
\noindent We can now state the next result. 
\begin{theorem}
\label{thm:markov_p}
\ifdefined\itistwocolumns
Suppose Assumptions~\ref{ass:feasibility}, \ref{ass:interpolation}, and \ref{ass:strong_convexity} hold. Furthermore, suppose \Cref{ass:ss_condition} with $\bar{\phi} \leq 1/L$. Let $x_1 \in \mathcal{X}$ and $\kappa = L/\mu$, then the iterates generated by eq.~\eqref{eq:sgd_ls} satisfy the following. For every $p \in (0,1]$, $n \geq 1$, and $\beta \in (0,1)$, with probability at least $1-\beta$, it holds that
\begin{align}\label{eq:hp_markov_p}
\Delta(x_{n+1}) &\leq \frac{1}{\beta^p} \cdot \frac{L \cdot \|x_{1} -x_{\star}\|^2}{2} \cdot \nonumber\\ &\bigg(1- a_p(\bar{\phi}, L) \cdot L \cdot \bar{\phi} \cdot \kappa^{-1} \bigg)^{pn},
\end{align}
\else
Suppose Assumptions~\ref{ass:feasibility}, \ref{ass:interpolation}, and \ref{ass:strong_convexity} hold. Furthermore, suppose \Cref{ass:ss_condition} with $\bar{\phi} \leq 1/L$. Let $x_1 \in \mathcal{X}$ and $\kappa = L/\mu$, then the iterates generated by eq.~\eqref{eq:sgd_ls} satisfy the following. For every $p \in (0,1]$, $n \geq 1$, and $\beta \in (0,1)$ it holds that
\begin{align}
\label{eq:hp_markov_p}
\bbP \bigg(\Delta(x_{n+1}) \leq \frac{1}{\beta^p} \cdot \frac{L \cdot \|x_{1} -x_{\star}\|^2}{2} \cdot \bigg(1- a_p(\bar{\phi}, L) \cdot \bar{\phi} \cdot L \cdot \kappa^{-1} \bigg)^{pn} \bigg) \geq 1-\beta,
\end{align}
\fi
where
\begin{align*}
a_p(\bar{\phi}, L) \coloneqq \frac{1-\big(1 - L \cdot \bar{\phi}\big)^{1/p}}{L \cdot \bar{\phi}}.
\end{align*}
\end{theorem}
The bound in eq.~\eqref{eq:hp_markov_p} exhibits a trade-off between the dependence on the failure probability $\beta$ and the number of iterations $n$. To put it into perspective, we note that the same analysis also yields a companion exponential rate in expectation, stated in Appendix~\ref{appendix:sc_ub_full}, of order $\exp\bigl(- \bar{\phi} \cdot L \cdot \kappa^{-1} \cdot n\bigr)$. The high-probability bound in \Cref{thm:markov_p} should be understood relative to this benchmark: for $p=1$, it preserves the in-expectation exponential rate in $n$ but incurs the poor confidence term $\beta^{-1}$, whereas for $p<1$ it improves the dependence on $\beta$ to $\beta^{-p}$ at the expense of a slower, yet still exponential, rate in $n$. For instance, if $\eta_k = 1/L$, then \Cref{prop:ss_condition_examples}~\cref{prop:ss_condition_eta} gives $\bar{\phi} = 1/L$, so the in-expectation rate is of order $\exp(- \kappa^{-1} \cdot n)$, while the corresponding high-probability rate becomes $\beta^{-p} \cdot \exp(- p \cdot \kappa^{-1} \cdot n)$. Thus, relative to the in-expectation benchmark, the confidence term improves as $p$ decreases, while the exponential rate constant deteriorates by a factor $p$. We also note that all the step-size policies in \Cref{prop:ss_condition_examples} satisfy $\bar{\phi} \leq 1/L$.

\noindent The next result shows that improving the dependence on the failure probability from $\beta^{-1}$ to $\beta^{-p}$ necessarily entails a deterioration of the exponential rate by a factor of order at least $p$.
\begin{theorem}
\label{thm:sc_lower_bound_p}
For every $p \in (0,1]$ and for every $\kappa > 1$, there exists a problem satisfying Assumptions~\ref{ass:feasibility}, \ref{ass:interpolation}, and \ref{ass:strong_convexity} with stochastic condition number $\kappa$ that satisfies the following. Let $x_1 \in \cX$ with $\Delta(x_1)>0$. Then for every step size policy and for every
\begin{align*}
a>-p\log\Bigl(1-\frac1\kappa\Bigr),
\end{align*}
there exist $n \ge 1$ and $\beta\in(0,1)$ s.t. for the iterates generated by eq.~\eqref{eq:sgd_ls}, with probability larger than $\beta$, one has
\begin{align}
\Delta(x_{n+1}) > \beta^{-p} \cdot \Delta(x_1) \cdot e^{-an}.
\end{align}
\end{theorem}
The above shows that the deterioration of the exponential rate in \Cref{thm:markov_p} is essentially tight. Indeed, the lower bound rules out any exponent larger than
\[
-p\log\Bigl(1-\frac1\kappa\Bigr) \leq \frac{p}{\kappa-1},
\]
which, for large $\kappa$, is of order at most $p/\kappa$. When $\eta_k=1/L$, \Cref{thm:markov_p} yields an upper bound with exponential rate of order $(p/\kappa)\cdot n$. Thus, in this case, the upper bound is asymptotically sharp as $\kappa\to\infty$. More generally, for any step-size policy such that $\bar{\phi} \in (0, 1/L)$, the exponent in \Cref{thm:markov_p} is of order $p/\kappa$, and is therefore optimal up to numerical constants.

\begin{remark}
\label{rmk:sc_lower_bound_p}
We make the following comments.
\begin{enumerate}[label={\rm(\roman*)}]
\item\label{rmk:sc_lower_bound_p_obst} 
An inspection of the proof in Appendix~\ref{subsec:appendix_sc_lower_bound_p} reveals that the obstruction comes from the presence of a ``lazy'' event, of measure exponentially small in $n$, on which the iterates do not move at all. Importantly, this mechanism is completely independent of the specific step-size policy: on that event, all sampled gradients vanish, so the iterate remains frozen regardless of how the step size is selected. Thus, the deterioration of the exponential rate exhibited by \Cref{thm:markov_p} is not a proof artifact, but a genuine obstruction. On the other hand, the worst-case example used in the proof still falls within the scope of \Cref{thm:markov_p}, and therefore SGD converges exponentially fast in expectation. 
\item\label{rmk:sc_lower_bound_p_polylog}
Importantly, eq.~\eqref{eq:hp_markov_p} is not obtained by a naive application of Markov's inequality. Rather, the key point is that under interpolation the moments of the error contract exponentially fast, which allows one to trade the dependence on $\beta$ against the exponential rate in $n$. It is therefore natural to ask whether one can retain an exponentially decaying dependence on $n$ while improving the confidence term to $\text{polylog(}\beta^{-1})$. The answer is negative in general. Indeed, for every $p>0$, every polylogarithmic factor in $\beta^{-1}$ is eventually dominated by $\beta^{-p}$ as $\beta\downarrow 0$. Since \Cref{thm:sc_lower_bound_p} holds for every $p\in(0,1]$, it follows that no bound with only polylogarithmic dependence on $1/\beta$ can hold simultaneously for every $n\geq 1$ and every $\beta\in(0,1)$, even for arbitrarily small exponential rate constant $a>0$.
\end{enumerate}
\end{remark}
\section{Extensions}\label{sec:extensions}
In this section, we present two extensions to the main results. Throughout this section, we only consider the case of constant step size $\eta_k = \eta$ for every $k \geq 1$, and leave for future work the generalization to the other policies. All proofs can be found in Appendix \ref{appendix:extensions}. 

\paragraph{Low Noise Regime.} In some applications, exact interpolation may fail, and the stochastic gradient noise may remain nonzero even at the optimum. Motivated by this, we consider the following relaxation of \Cref{ass:interpolation}.
\begin{assumption}[Low Noise Regime]\label{ass:exp_noise_optimum}
There exists $x_{\star} \in \cX_{\star}$ and $\sigma_{\star} > 0$ s.t.
\begin{align*}
\bbE\left[ \exp\left( \frac{\|\nabla f(x_\star,Z)\|^2}{\sigma_\star^2} \right) \right] \le 2.
\end{align*}
\end{assumption}
Previous works focusing on in-expectation guarantees typically assume instead that the noise at the optimum has finite variance \cite{Attia2025,Garrigos2025}. Here we strengthen this requirement by imposing a subGaussian condition at $x_\star$. This assumption is still local, in that it only concerns the noise at the optimum, yet it is strong enough to yield the following high-probability result. Therefore, it should not be confused with the uniform subGaussian noise conditions commonly used in high-probability analyses of SGD \cite{Liu2022,Eldowa2024,Liu2024}. Those works require a tail bound on $\|\nabla f(x,Z)-\nabla F(x)\|$ with a parameter that is uniform over all $x \in \cX$. In this sense, the resulting setting is more general, and correspondingly harder to analyze.

\begin{theorem}\label{thm:hp_generalized_interpolation}
Suppose Assumptions \ref{ass:feasibility}, \ref{ass:convexity}, and \ref{ass:exp_noise_optimum} hold. Let $x_1 \in \cX$, and for every $n \geq 1$ let
\begin{align*}
\eta
=
\min\left\{
\frac{1}{4 \cdot L},
\frac{\|x_1-x_\star\|}{32 \cdot \sigma_\star \cdot \sqrt{n \cdot \log(3/\beta)}}
\right\}.
\end{align*}
Then the iterates generated by eq. \eqref{eq:sgd_ls} satisfy the following. For $\beta \in (0,1)$, with probability at least $1-\beta$, it holds that
\ifdefined\itistwocolumns
\begin{align*}
\Delta(\bar x_n) &\le 8 \cdot \bigg(2 + \log \Big(\frac 3 \beta\Big)\bigg) \cdot \frac{L \cdot \|x_1-x_\star\|^2}{n}\\
&\quad+ 98 \cdot  \sigma_\star \cdot \|x_1-x_\star\| \cdot 
\sqrt{\frac{\log\Big( \frac 3 \beta \Big)}{n}}.
\end{align*}
\else
\begin{align*}
\Delta(\bar x_n) \le 8 \cdot \bigg(2 + \log \Big(\frac 3 \beta\Big)\bigg) \cdot \frac{L \cdot \|x_1-x_\star\|^2}{n}
+ 98 \cdot  \sigma_\star \cdot \|x_1-x_\star\| \cdot 
\sqrt{\frac{\log\Big( \frac 3 \beta \Big)}{n}}.
\end{align*}
\fi
\end{theorem}
We note that the above result is useful in two respects: it recovers the same order as \Cref{thm:main_result} in the limit $\sigma_{\star} \to 0$, and it matches the best known in-expectation guarantees up to polylogarithmic factors in $1/\beta$. Moreover, we stress that our choice of $\eta$ is not an artifact of the high-probability analysis: existing in-expectation results in the low-noise regime also require tuning in terms of $n$, $\sigma_\star$, $L$, and $\|x_1-x_\star\|$ \cite{Attia2025,Garrigos2025}. An interesting direction for future work is to design adaptive methods that avoid this prior knowledge. At the proof level, the result relies on ideas that go beyond those used for \Cref{thm:main_result}. The first difficulty is that, in the low-noise regime, SGD no longer enjoys Fej\'er monotonicity, so one must explicitly control the location of the iterates. The second is that $\Delta(x_k,Z_k)$ does not vanish even when $x_k=x_\star$, and therefore its fluctuations cannot be treated directly as in the interpolation regime. To overcome this, we introduce a first-order correction term that restores an interpolation-type structure and yields a fluctuation term whose variance can still be absorbed into the optimization error. This is what gives the $n^{-1}$ contribution. The genuinely new term comes from the residual noise at the optimum, which induces an additional martingale of size $\sigma_\star \cdot R_\eta \cdot \sqrt{\log(1/\beta)/n}$, where $R_\eta$ is a high-probability bound on the radius of the iterate trajectory. Choosing $\eta$ so that $R_\eta$ remains of order $\|x_1-x_\star\|$ then yields the two-regime bound, consisting of an interpolation-like term of order $n^{-1}$ and a noise-driven term of order $\sigma_\star n^{-1/2}$. The choice of this well-tuned step size is a consequence of a more general bound that would allow $\eta\in (0,2/L)$; details are deferred to the appendix.

\paragraph{Nonconvex Optimization.} A broad class of modern optimization problems arising in machine learning is nonconvex. We therefore consider the following relaxation of \Cref{ass:strong_convexity}.  
\begin{assumption}
\label{ass:nonconvexity}
For $\cP$-almost every $z \in \cZ$, $f(\cdot, z)\colon \mathcal{X} \rightarrow \mathbb{R}$ is $L$-smooth (this implies that $F$ is $L$-smooth). Moreover, there exists $\mu > 0$ s.t.:
\begin{enumerate}[label={\rm(\roman*)}]
\item\label{ass:nonconvexity_PL} $\forall x \in \cX, \|\nabla F(x) \|^2 \geq 2 \cdot \mu \cdot \Delta(x)$,
\item\label{ass:nonconvexity_common_minima} For $\cP$-almost every $z$, $\cX_{\star} \subseteq \arg\min_{x \in \cX} f(x, z)$.
\end{enumerate}
\end{assumption}
\noindent Assumption~\ref{ass:nonconvexity}\ref{ass:nonconvexity_PL} is the classical Polyak--\L ojasiewicz condition \cite{Polyak1963,Lojasiewicz1963}, which does not require convexity of $F$, and is satisfied by several genuinely nonconvex optimization problems of interest in machine learning \cite{Liu2022}. Assumption~\ref{ass:nonconvexity}\ref{ass:nonconvexity_common_minima} strengthens interpolation (\Cref{ass:interpolation}), and requires that every global minimizer of $F$ also minimizes each sample loss. This is natural, for instance, in interpolation-type regimes where all sample losses are non-negative and can be driven to zero simultaneously. We can now state our guarantees for SGD under Assumption~\ref{ass:nonconvexity}. The original statement is provided in \Cref{thm:pl_qg_extension_full}, while we report here a simplified version obtained by fixing $\eta = (1/5) \cdot (\mu/L^2)$.

\begin{theorem}
\label{thm:pl_qg_extension}
Suppose Assumptions~\ref{ass:feasibility} and \ref{ass:nonconvexity} hold. Let $x_1 \in \cX$ and $\eta = (1/5) \cdot (\mu/L^2)$. Then the iterates generated by eq.~\eqref{eq:sgd_ls} satisfy the following. For every $n\geq 1$ and every $\beta\in(0,1)$, with probability at least $1-\beta$, it holds that
\ifdefined\itistwocolumns
\begin{equation}\label{eq:pl_qg_hp_bound_main_spec}
    \begin{split}
\Delta(x_{n+1}) \leq \frac{1}{\sqrt{\beta}} \cdot \Delta(x_1) \cdot \left( 1 - \frac 1 5 \cdot \left(\frac{\mu}{L}\right)^2 \right)^{n/2}.
    \end{split}
\end{equation}
\else
\begin{align}
\label{eq:pl_qg_hp_bound_main_spec}
\bbP\bigg( \Delta(x_{n+1}) \leq \frac{1}{\sqrt{\beta}} \cdot \Delta(x_1) \cdot ( 1 - (1/5) \cdot (\mu/L)^2 )^{n/2} \bigg) \geq 1-\beta.
\end{align}
\fi
\end{theorem}
The above result improves upon the application of Markov inequality to the in-expectation rate (see \Cref{appendix:extensions_nonconvex_optimization})
in terms of the dependence on $\beta$, and implies a rate of order $\beta^{-1/2} \cdot \exp(-(1/10) \cdot \kappa^{-2} \cdot n)$, for $\kappa = L/\mu$. We note that, even if smooth strongly convex functions are a subset of functions satisfying \cref{ass:nonconvexity}, \Cref{thm:sc_lower_bound_p} does not yet clarify whether this rate is tight, since it features a worse dependence on the condition number. 

\section{Conclusions}\label{sec:conclusions}
Under the interpolation assumption, we establish high-probability fast rates for SGD that match the corresponding in-expectation guarantees, without requiring additional assumptions, for $L$-smooth convex and strongly convex functions. Beyond interpolation, in the convex low-noise regime, we recover the expected $n^{-1}+\sigma_\star \cdot n^{-1/2}$ behavior under the sole assumption that the noise at the solution is $\sigma_\star$-subGaussian. Going beyond convexity, we also obtain an exponential high-probability rate for $\mu$-PL functions (with samplewise common minimizers), analogous to the one proved in the strongly convex case. 

In both of these latter settings, however, the resulting step-sizes are not adaptive, as they depend on unknown quantities such as $\mu$, $L$, and the time horizon $n$. Developing adaptive procedures that remove this prior knowledge is therefore an important direction for future work. Other directions that we plan to explore include last-iterate convergence, and extensions to mirror descent. 

\ifdefined\itistwocolumns
\subsection*{AI use statement}
In this work, we used generative AI tools for editing purposes. We have not used generative AI tools to prove or derive the mathematical statements. We take responsibility for the final content of this work, including text, claims, or artifacts produced with the aid of generative AI.
\else
\subsection*{AI use statement}
In this work, we used generative AI tools exclusively for editing purposes, more specifically for spelling and grammar corrections throughout the manuscript and for rephrasing selected passages in Section 1. We have not used generative AI tools to prove or derive the mathematical statements. Moreover, we have not used AI tools to present, interpret and comment our mathematical results. 
We take responsibility for the final content of this work, including text, claims, or artifacts produced with the aid of generative AI.
\fi

\bibliographystyle{alpha}
\bibliography{refs}

@article{Apidopoulos2022,
  title={Convergence rates for the heavy-ball continuous dynamics for non-convex optimization, under Polyak--{\L}ojasiewicz condition},
  author={Apidopoulos, Vassilis and Ginatta, Nicolo and Villa, Silvia},
  journal={Journal of Global Optimization},
  volume={84},
  number={3},
  pages={563--589},
  year={2022},
  publisher={Springer}
}

@inproceedings{Attia2025,
  title={Fast Last-Iterate Convergence of {SGD} in the Smooth Interpolation Regime},
  author={Attia, Amit and Schliserman, Matan and Sherman, Uri and Koren, Tomer},
  booktitle={The Thirty-ninth Annual Conference on Neural Information Processing Systems},
  year={2025},
}

@article{Agmon1954,
  title={The relaxation method for linear inequalities},
  author={Agmon, Shmuel},
  journal={Canadian Journal of Mathematics},
  volume={6},
  pages={382--392},
  year={1954},
  publisher={Cambridge University Press}
}

@article{Agarwal2012,
    author={Agarwal, Alekh and Bartlett, Peter L. and Ravikumar, Pradeep and Wainwright, Martin J.},
    title={Information-Theoretic Lower Bounds on the Oracle Complexity of Stochastic Convex Optimization}, 
    journal={IEEE Transactions on Information Theory}, 
    year={2012},
    volume={58},
    number={5},
    pages={3235-3249},
}

@book{Bauschke2011,
  title={Convex Analysis and Monotone Operator Theory in Hilbert Spaces},
  author={Bauschke, Heinz H and Combettes, Patrick L},
  publisher={Springer},
  year={2011}
}

@article{Bernstein1924,
    author={Bernstein, S. N.},
    title={On a modification of Chebyshev's inequality and of the error formula of Laplace (in Russian)},
    journal={Uch. zap. n.-i., kafedr Ukraini, Otd. matem},
    pages={38--49},
    year={1924}
}

@article{Blumer1989,
  title={Learnability and the Vapnik-Chervonenkis dimension},
  author={Blumer, Anselm and Ehrenfeucht, Andrzej and Haussler, David and Warmuth, Manfred K},
  journal={Journal of the ACM (JACM)},
  volume={36},
  number={4},
  pages={929--965},
  year={1989},
  publisher={ACM New York, NY, USA}
}

@article{Bottou2018,
  title={Optimization methods for large-scale machine learning},
  author={Bottou, L{\'e}on and Curtis, Frank E and Nocedal, Jorge},
  journal={SIAM review},
  volume={60},
  number={2},
  pages={223--311},
  year={2018},
  publisher={SIAM}
}

@article{Doob1940,
  title={Regularity properties of certain families of chance variables},
  author={Doob, Joseph L},
  journal={Transactions of the American Mathematical Society},
  volume={47},
  number={3},
  pages={455--486},
  year={1940}
}

@inproceedings{Eldowa2024,
  title={General tail bounds for non-smooth stochastic mirror descent},
  author={Eldowa, Khaled and Paudice, Andrea},
  booktitle={International Conference on Artificial Intelligence and Statistics},
  pages={3205--3213},
  year={2024},
  organization={PMLR}
}

@inproceedings{Karimi2016,
  title={Linear convergence of gradient and proximal-gradient methods under the polyak-{\l}ojasiewicz condition},
  author={Karimi, Hamed and Nutini, Julie and Schmidt, Mark},
  booktitle={Joint European conference on machine learning and knowledge discovery in databases},
  pages={795--811},
  year={2016},
  organization={Springer}
}

@incollection{Fejer1922,
  title={{\"U}ber die Lage der Nullstellen von Polynomen, die aus Minimumforderungen gewisser Art entspringen},
  author={Fej{\'e}r, Leopold},
  booktitle={Mathematische Annalen},
  pages={41--48},
  volume={85},
  issue={1},
  year={1922},
  publisher={Springer}
}

@article{Galli2023,
  title={Don't be so monotone: Relaxing stochastic line search in over-parameterized models},
  author={Galli, Leonardo and Rauhut, Holger and Schmidt, Mark},
  journal={Advances in Neural Information Processing Systems},
  volume={36},
  pages={34752--34764},
  year={2023}
}

@article{Garrigos2025,
  title={Last-iterate complexity of {SGD} for convex and smooth stochastic problems},
  author={Garrigos, Guillaume and Cortild, Daniel and Ketels, Lucas and Peypouquet, Juan},
  journal={arXiv preprint arXiv:2507.14122},
  year={2025}
}

@inproceedings{Gorbunov2020,
    title={Stochastic optimization with heavy-tailed noise via accelerated gradient clipping},
    author={Gorbunov, Eduard and Danilova, Marina and Gasnikov, Alexander},
    booktitle={Proceedings of the 34th Annual Conference on Neural Information Processing Systems},
    pages={15042--15053},
    year={2020}
}

@inproceedings{Harvey2019,
  title={Tight analyses for non-smooth stochastic gradient descent},
  author={Harvey, Nicholas JA and Liaw, Christopher and Plan, Yaniv and Randhawa, Sikander},
  booktitle={Conference on Learning Theory},
  pages={1579--1613},
  year={2019},
  organization={PMLR}
}

@article{Jain2021,
    title = {Making the Last Iterate of {SGD} Information Theoretically Optimal},
    author = {Jain, Prateek and Nagaraj, Dheeraj M. and Netrapalli, Praneeth},
    journal = {SIAM Journal on Optimization},
    volume = {31},
    number = {2},
    pages = {1108-1130},
    year = {2021},
}

@article{Jin2019,
  title={A short note on concentration inequalities for random vectors with subgaussian norm},
  author={Jin, Chi and Netrapalli, Praneeth and Ge, Rong and Kakade, Sham M and Jordan, Michael I},
  journal={arXiv preprint arXiv:1902.03736},
  year={2019}
}

@article{Jin2024,
  title={High probability complexity bounds for adaptive step search based on stochastic oracles},
  author={Jin, Billy and Scheinberg, Katya and Xie, Miaolan},
  journal={SIAM Journal on Optimization},
  volume={34},
  number={3},
  pages={2411--2439},
  year={2024},
  publisher={SIAM}
}

@article{Klochkov2021,
  title={Stability and deviation optimal risk bounds with convergence rate $ O (1/n) $},
  author={Klochkov, Yegor and Zhivotovskiy, Nikita},
  journal={Advances in Neural Information Processing Systems},
  volume={34},
  pages={5065--5076},
  year={2021}
}

@article{Liu2022,
	title={Loss landscapes and optimization in over-parameterized non-linear systems and neural networks},
	author={Liu, Chaoyue and Zhu, Libin and Belkin, Mikhail},
	journal={Applied and Computational Harmonic Analysis},
	volume={59},
	pages={85--116},
	year={2022},
	publisher={Elsevier}
}

@inproceedings{Liu2023,
  title={High probability convergence of stochastic gradient methods},
  author={Liu, Zijian and Nguyen, Ta Duy and Nguyen, Thien Hang and Ene, Alina and Nguyen, Huy},
  booktitle={International conference on machine learning},
  pages={21884--21914},
  year={2023},
  organization={PMLR}
}

@inproceedings{Liu2024,
    author={Liu, Zijian and Zhou, Zhengyuan},
    title={Revisiting the last-iterate convergence of stochastic gradient methods},
    booktitle={Proceedings of the 12-th International Conference on Learning Representations (to appear)},
    year={2024}
}

@inproceedings{Loizou2021,
  title={Stochastic polyak step-size for {SGD}: An adaptive learning rate for fast convergence},
  author={Loizou, Nicolas and Vaswani, Sharan and Laradji, Issam Hadj and Lacoste-Julien, Simon},
  booktitle={International conference on artificial intelligence and statistics},
  pages={1306--1314},
  year={2021},
  organization={PMLR}
}

@article{Lojasiewicz1963,
    author={{{{{\L}ojasiewicz, S.}}}},
    title={A topological property of real analytic subsets (in French)},
    journal={Coll. du CNRS, Les équations},
    pages={87--89},
    year={1963},
}

@article{Lou2022,
    title={Beyond sub-gaussian noises: Sharp concentration analysis for stochastic gradient descent},
    author={Lou, Zhipeng and Zhu, Wanrong and Wu, Wei Biao},
    journal={The Journal of Machine Learning Research},
    volume={23},
    number={1},
    pages={2227-2248},
    year={2022}
}

@inproceedings{Ma2018,
  title={The power of interpolation: Understanding the effectiveness of {SGD} in modern over-parametrized learning},
  author={Ma, Siyuan and Bassily, Raef and Belkin, Mikhail},
  booktitle={International Conference on Machine Learning},
  pages={3325--3334},
  year={2018},
  organization={PMLR}
}

@article{Madden2024,
  title={High probability convergence bounds for non-convex stochastic gradient descent with sub-weibull noise},
  author={Madden, Liam and Dall'Anese, Emiliano and Becker, Stephen},
  journal={Journal of Machine Learning Research},
  volume={25},
  number={241},
  pages={1--36},
  year={2024}
}

@article{Moulines2011,
  title={Non-asymptotic analysis of stochastic approximation algorithms for machine learning},
  author={Moulines, Eric and Bach, Francis},
  journal={Advances in neural information processing systems},
  volume={24},
  year={2011}
}

@article{Bach2013,
  title={Non-strongly-convex smooth stochastic approximation with convergence rate O (1/n)},
  author={Bach, Francis and Moulines, Eric},
  journal={Advances in neural information processing systems},
  volume={26},
  year={2013}
}

@book{Nemirovski1983,
    title={Problem complexity and method efficiency in optimization},
    author={Nemirovskij, Arkadij Semenovi{\v{c}} and Yudin, David Borisovich},
    year={1983},
    publisher={Wiley-Interscience}
}

@article{Nemirovski2009,
    author = {Nemirovski, A. and Juditsky, A. and Lan, G. and Shapiro, A.},
    title = {Robust Stochastic Approximation Approach to Stochastic Programming},
    journal = {SIAM Journal on Optimization},
    volume = {19},
    number = {4},
    pages = {1574-1609},
    year = {2009}
}

@book{Nesterov2004,
	title={Introductory lectures on convex optimization: A basic course},
	author={Nesterov, Yurii},
	year={2004},
	publisher={Appl. Optim. 87,	Kluwer Academic Publishers, Dordrecht, The Netherlands}
}

@inproceedings{Nguyen2023,
    author={Nguyen, Ta Duy and Ene, Alina and Nguyen, Huy L.},
    title={Improved Convergence in High Probability of Clipped Gradient Methods with Heavy Tails},
    booktitle={Proceedings of the 37th Annual Conference on Neural Information Processing Systems (to appear)},
    year={2023}
}

@article{Parletta2024,
    author={D. A. Parletta and A. Paudice and M. Pontil and S. Salzo},
    title={High probability bounds for stochastic subgradient schemes with heavy-tailed noise},
    journal={SIAM Journal on Mathematics of Data Science (to appear)},
    year={2024}
}

@article{Polyak1963,
  title={Gradient methods for minimizing functionals (in Russian)},
  author={Polyak, Boris Teodorovich},
  journal={Zhurnal vychislitel'noi matematiki i matematicheskoi fiziki},
  volume={3},
  number={4},
  pages={643--653},
  year={1963},
  publisher={Russian Academy of Sciences, Branch of Mathematical Sciences}
}

@article{Robbins1951,
  title={A stochastic approximation method},
  author={Robbins, Herbert and Monro, Sutton},
  journal={The annals of mathematical statistics},
  pages={400--407},
  year={1951},
  publisher={JSTOR}
}

@inproceedings{Sebbouh2021,
  title={Almost sure convergence rates for stochastic gradient descent and stochastic heavy ball},
  author={Sebbouh, Othmane and Gower, Robert M and Defazio, Aaron},
  booktitle={Conference on Learning Theory},
  pages={3935--3971},
  year={2021},
  organization={PMLR}
}

@inproceedings{Vaswani2019a,
  title={Fast and faster convergence of {SGD} for over-parameterized models and an accelerated perceptron},
  author={Vaswani, Sharan and Bach, Francis and Schmidt, Mark},
  booktitle={The 22nd international conference on artificial intelligence and statistics},
  pages={1195--1204},
  year={2019},
  organization={PMLR}
}

@article{Vaswani2019b,
  title={Painless stochastic gradient: Interpolation, line-search, and convergence rates},
  author={Vaswani, Sharan and Mishkin, Aaron and Laradji, Issam and Schmidt, Mark and Gidel, Gauthier and Lacoste-Julien, Simon},
  journal={Advances in neural information processing systems},
  volume={32},
  year={2019}
}

@book{Vershynin2018,
  title={High-dimensional probability: An introduction with applications in data science},
  author={Vershynin, Roman},
  volume={47},
  year={2018},
  publisher={Cambridge university press}
}

@article{Young1912,
  title={On classes of summable functions and their Fourier series},
  author={Young, William Henry},
  journal={Proceedings of the Royal Society of London A},
  volume={87},
  number={594},
  pages={225--229},
  year={1912},
  publisher={The Royal Society London}
}

\makeatletter
\let\etalchar\relax
\makeatother


\newpage
\appendix
\ifdefined\itistwocolumns

\else
\newpage
\appendix
\fi
\section{Appendix: Mathematical Facts and Technical Lemmas}\label{appendix:mathematical_facts}
\paragraph{Random Variables.} We recall the definition of subGaussian random variable \cite{Vershynin2018}.
\begin{definition}(SubGaussian Random Variables)
A real valued random variable $X$ is subGaussian with parameter $\sigma > 0$ iff
\[
\bbE[\exp(X^2/\sigma^2)] \leq 2.
\]
\end{definition}
By extension to vectors we obtain the following.
\begin{definition}(SubGaussian Random Vectors)
Let $\cX$ be a real separable Hilbert Space. An $\cX$-valued random variable $X$ is:
\begin{itemize}
\item subGaussian with parameter $\sigma > 0$ iff for every $u \in \cX$, $\langle X, u \rangle$ is $\sigma \cdot \|u\|$-\text{subGaussian}.

\item norm subGaussian with parameter $\sigma > 0$ iff $\|X\|$ is $\sigma$-subGaussian.
\end{itemize}
\end{definition}
The next lemma shows that a norm sub-Gaussian random vector is also sub-Gaussian along every linear form. This fact can be found, for instance, in \cite[Lemma 3]{Jin2019}; we include a proof here in order to keep track of explicit numerical constants.
\begin{lemma}\label{lem:exp_square_implies_linear_subgaussian}
Let $\cX$ be a real Hilbert space, and let $X$ be a $\cX$-valued random variable. Suppose that $X$ is mean-zero and norm sub-Gaussian with parameter $\sigma > 0$. Then, for every $u\in\cX$ and every $\lambda\in\bbR$,
\begin{align*}
\bbE\left[ \exp\!\big(\lambda \langle X,u\rangle\big)\right] \leq \exp\left( 2\log(2)\,\lambda^2\sigma^2\|u\|^2\right).
\end{align*}
\end{lemma}

\begin{proof}
The claim is trivial if $u=0$, so assume $u\neq 0$, and define $Y \coloneqq \langle X,u\rangle$. Then $\bbE[Y]=\langle \bbE[X],u\rangle=0$, and by Cauchy--Schwarz, $|Y|\le \|u\|\,\|X\|$.
Hence
\begin{align}\label{eq:lemma_norm_product_1}
\bbE\!\left[
\exp\!\left(
\frac{Y^2}{\sigma^2\|u\|^2}
\right)
\right]
\le
\bbE\!\left[
\exp\!\left(
\frac{\|X\|^2}{\sigma^2}
\right)
\right]
\le 2.
\end{align}

We now distinguish two cases.

\medskip
\noindent\emph{Case 1:} $|\lambda|\le (\sigma\|u\|)^{-1}$.

Let $Y'$ be an independent copy of $Y$. Then
\begin{align*}
\bbE[e^{\lambda Y}] &= \bbE\left[\exp(\lambda(Y-\bbE[Y']))\right] \\
&\overset{(a)}{\leq} \bbE\left[\exp(\lambda(Y-Y'))\right] \\
&\overset{(b)}{=} \bbE\left[\cosh(\lambda(Y-Y'))\right] \\
&\overset{(c)}{\leq} \bbE\left[\exp\left(\frac{\lambda^2(Y-Y')^2}{2} \right) \right] \\
&\overset{(d)}{\leq} \bbE[\exp(\lambda^2Y^2)] \bbE[\exp(\lambda^2(Y')^2)] \\
&= \left(\bbE[e^{\lambda^2Y^2}]\right)^2 \\
&= \Big(\bbE\left[ \left( \exp(Y^2/(\sigma^2\|u\|^2)) \right)^{\lambda^2\sigma^2\|u\|^2} \right]\Big)^2 \\
&\overset{(e)}{\leq} \Big(\left( \bbE\left[ \exp(Y^2/(\sigma^2\|u\|^2)) \right] \right)^{\lambda^2\sigma^2\|u\|^2}\Big)^2 \\
&\overset{(f)}{\leq} 2^{2\lambda^2\sigma^2\|u\|^2} \\
&= \exp\left( 2\log(2)\,\lambda^2\sigma^2\|u\|^2 \right),
\end{align*}
where: (a) follows from Jensen's inequality, (b) follows from the symmetry of $Y-Y'$, (c) follows from the elementary inequality $\cosh(t)\le e^{t^2/2}$, (d) follows by $(a-b)^2 \leq 2 (a^2 + b^2)$, (e) follows from the concavity of $t\mapsto t^{\lambda^2\sigma^2\|u\|^2}$ on $(0,\infty)$, and (f) follows by eq.~\eqref{eq:lemma_norm_product_1}.

\medskip
\noindent\emph{Case 2:} $|\lambda|>(\sigma\|u\|)^{-1}$.

For every $y\in\bbR$,
\[
\lambda y \leq \frac{\lambda^2\sigma^2\|u\|^2}{2} + \frac{y^2}{2\sigma^2\|u\|^2}.
\]
Therefore
\[
e^{\lambda Y} \leq \exp\left( \frac{\lambda^2\sigma^2\|u\|^2}{2} \right) \exp\left( \frac{Y^2}{2\sigma^2\|u\|^2} \right).
\]
Taking expectations and using Cauchy--Schwarz gives
\[
\bbE[e^{\lambda Y}] \leq \exp\left( \frac{\lambda^2\sigma^2\|u\|^2}{2} \right) \left( \bbE\left[ \exp\!\left( \frac{Y^2}{\sigma^2\|u\|^2} \right) \right] \right)^{1/2} \leq
\sqrt{2} \exp\!\left( \frac{\lambda^2\sigma^2\|u\|^2}{2} \right).
\]
Since $|\lambda|>(\sigma\|u\|)^{-1}$, we have
\[
\frac{1}{2}\log(2) \le \frac{\lambda^2\sigma^2\|u\|^2}{2}\log(2),
\]
and therefore
\[
\sqrt{2} \exp\left( \frac{\lambda^2\sigma^2\|u\|^2}{2} \right) \leq \exp\left( \frac{1+\log(2)}{2}\lambda^2\sigma^2\|u\|^2 \right) \leq \exp\left( 2\log(2)\,\lambda^2\sigma^2\|u\|^2 \right),
\]
because
\[
\frac{1+\log(2)}{2}\le 2\log(2).
\]
\noindent Combining the two cases concludes the proof.
\end{proof}

The following is a classical result \cite{Doob1940}.
\begin{lemma}[Doob's Maximal Inequality]
\label{lem:doob_maximal}
Let $(Y_m)_{m=0}^n$ be a nonnegative supermartingale. Then, for every $t>0$,
\begin{align*}
\bbP\left(\max_{0 \leq m \leq n } Y_m \ge t\right) \leq \frac{\bbE[Y_0]}{t}.    
\end{align*}
\end{lemma}

\paragraph{Convex Analysis.} The following two lemmas are often used in our results, we report here the statement, while for the proofs we reference for instance to \cite{Bauschke2011}.

\begin{lemma}[Descent Lemma] Suppose that $F:\cX \to \R$ is $L$-smooth. Then for any $x,y \in \cX$, we have
	\begin{equation*}
		F(y) \leq F(x) + \langle\nabla F(x), y-x\rangle + \frac{L}{2} \|x-y\|^2.
	\end{equation*}
\end{lemma}

\begin{lemma}[Co-coercivity]
	Suppose that $F:\cX \to \R$ is $L$-smooth and convex, then $\forall x, y \in \cX$, we have both
	\begin{align*}
		\frac{1}{2L} ||\nabla F(y)-\nabla F(x)||^2 &\leq F(y)-F(x) -\langle\nabla F(x), y-x\rangle\\
		\frac{1}{L} ||\nabla F(y)-\nabla F(x)||^2&\leq \langle\nabla F(y)-\nabla F(x),y-x\rangle
	\end{align*}
\end{lemma}

\newpage
\section{Appendix: Section \ref{sec:convex_smooth_optimization}}\label{sec:appendix_convex_smooth_optimization}
\subsection{Stochastic Line Search}\label{subsec:appendix_sls}
\begin{algorithm}[t]
\caption{Monotone Stochastic Line Search}
\label{alg:backtrack}
\begin{algorithmic}[1]
\STATE{\bfseries Input:} $x_k \in \mathcal{X}, Z_k \sim \cP, c \in (1/2, 1), \delta \in (0,1), \eta > 0$ 
\STATE $l_k = 0$
\STATE $\eta_k = \eta$
\WHILE{$
f(x_{k}-\eta_k \cdot \nabla f(x_k,Z_k), Z_k) > f(x_k,Z_k) - c \cdot \eta_k \cdot \|\nabla f(x_k, Z_k)\|^2$}  
\STATE $l_k = l_k +1$
\STATE $\eta_k = \delta^{l_k} \cdot \eta $
\ENDWHILE
\STATE{\bfseries Return:} $\eta_k$
\end{algorithmic}
\end{algorithm}
\noindent In this section we analyze the line search procedure described in \Cref{alg:backtrack}. The following already appeared for instance in Lemma 1 of \cite{Vaswani2019b} and we report it for completeness.
\begin{lemma}
\label{lemma:sls}
Suppose \Cref{ass:convexity} holds. Then, for every $k \geq 1$, the output of \Cref{alg:backtrack} satisfies
\begin{align}
\label{eq:ss_lb}
\forall k \in \mathbb{N}: \quad \eta_k \geq \min\bigg\{\frac{2 \cdot (1-c) \cdot \delta}{L}, \eta \bigg\} \eqqcolon \bar{\eta},
\end{align}
and note that, since $ c < 1$ and $\eta > 0$, $\eta_k \geq \bar{\eta} > 0$. 
\end{lemma}
\begin{proof}
We start noting that, for every $k \geq 1$ and for every $\eta>0$, we have
\begin{align}
\label{eq:ss_lb_analysis}
f (x_k -\eta \cdot \nabla f(x_k, Z_k), Z_k) &\overset{(a)}{\leq} f(x_k, Z_k) + \langle \nabla f(x_k, Z_k), x_k -\eta \cdot \nabla f(x_k, Z_k) - x_k \rangle \nonumber \\
&\quad + \frac{\eta^2L}{2} \| \nabla f(x_k, Z_k) \|^2 \nonumber \\
&= f(x_k, Z_k) - \left(\eta - \frac{\eta^2 \cdot L}{2} \right) \|\nabla f(x_k, Z_k)\|^2,
\end{align}
where (a) follows from the smoothness of $f(\cdot, Z_k)$ (\Cref{ass:convexity}). For $\eta > 0$, define 
\begin{align*}
p_k(\eta) &\coloneqq f(x_k, Z_k) - c \cdot \eta \cdot \| \nabla f(x_k, Z_k) \|^2 \\
q_k(\eta) &\coloneqq f(x_k, Z_k) - \left(\eta - \frac{\eta^2 \cdot L}{2} \right) \| \nabla f(x_k, Z_k) \|^2,
\end{align*}
and note that eq.~\eqref{eq:ss_lb_analysis} can be rewritten as
\begin{align}
\label{eq:q_eta}
f (x_k -\eta \cdot \nabla f(x_k, Z_k), Z_k) \leq q_k(\eta).
\end{align}
Now, \Cref{alg:backtrack} admits two possible outcomes. If $l_k=0$, then $\eta_k = \eta$. Otherwise, $l_k>0$. In this case, we have $\eta_k < \eta$. Moreover, by the condition in line 4 of \Cref{alg:backtrack}, we have 
\begin{align}
\label{eq:p_eta}
f(x_k - (\eta_k/\delta) \nabla f(x_k, Z_k), Z_k) &> f(x_k, Z_k) - c \cdot \frac{\eta_k}{\delta} \cdot  \|\nabla f(x_k, Z_k)\|^2  \nonumber \\
&= p_k(\eta_k/\delta).
\end{align}
Combining eq.~\eqref{eq:q_eta} with  eq.~\eqref{eq:p_eta} yields $p_k(\frac{\eta_k}{\delta})\leq q_k(\frac{\eta_k}{\delta})$. This implies that
\begin{equation*}
f(x_k, Z_k) - c \cdot \frac{\eta_k}{\delta} \cdot \| \nabla f(x_k, Z_k) \|^2 \leq f(x_k, Z_k) - \left(\frac{\eta_k}{\delta} - \frac{\eta_k^2 \cdot L}{2 \cdot \delta^2} \right) \cdot \| \nabla f(x_k, Z_k) \|^2
\end{equation*}
and consequently 
\begin{equation*}
- c \leq - \left(1 - \frac{\eta_k \cdot L}{2 \cdot \delta} \right) \Leftrightarrow \eta_k \geq \frac{2 \cdot \delta \cdot (1-c)}{L},
\end{equation*}
which proves eq.~\eqref{eq:ss_lb}.
\end{proof}
\begin{remark}
\label{rmk:termination}
\Cref{alg:backtrack} always halts in a finite number of steps. Indeed, by eq.~\eqref{eq:ss_lb_analysis}, any trial step size $\tilde{\eta}>0$ s.t.
\[
\tilde{\eta} \le \frac{2 \cdot (1-c)}{L}
\]
satisfies the condition in the while loop of \Cref{alg:backtrack}. Since the backtracking procedure tests the geometrically decreasing sequence
\[
\eta,\ \delta\eta,\ \delta^2\eta,\dots,
\]
with $\delta \in (0,1)$, it eventually produces a trial step size not larger than $2(1-c)/L$. Hence the while loop terminates after finitely many iterations.
\end{remark}

\subsection{Proof of Proposition \ref{prop:ss_condition_examples}}\label{subsec:appendix_convex_smooth_optimization_prop_ss_condition_examples}
Fix any $k\geq 1$. By the update rule in eq.~\eqref{eq:sgd_ls}, one has
\begin{align}
\label{eq:appendix_ss_general_expansion}
\|x_{k+1}-x_\star\|^2
&=
\|x_k-x_\star\|^2
-2\eta_k\langle \nabla f(x_k,Z_k),x_k-x_\star\rangle
+\eta_k^2\|\nabla f(x_k,Z_k)\|^2.
\end{align}
Since $f(\cdot,Z_k)$ is convex and $L$-smooth, and $x_\star\in \arg\min f(\cdot,Z_k)$ by interpolation, the standard cocoercivity inequality yields
\begin{align}
\label{eq:appendix_ss_cocoercive}
\Delta(x_k,Z_k)+\frac{1}{2L}\|\nabla f(x_k,Z_k)\|^2
\le
\langle \nabla f(x_k,Z_k),x_k-x_\star\rangle.
\end{align}
Combining eq.~\eqref{eq:appendix_ss_general_expansion} with eq.~\eqref{eq:appendix_ss_cocoercive}, we obtain
\begin{align}
\label{eq:appendix_ss_master}
\|x_{k+1}-x_\star\|^2
\le
\|x_k-x_\star\|^2
-2\eta_k\Delta(x_k,Z_k)
+\eta_k\Bigl(\eta_k-\frac{1}{L}\Bigr)\|\nabla f(x_k,Z_k)\|^2.
\end{align}
Moreover, by $L$-smoothness of $f(\cdot,Z_k)$ and the fact that $x_\star$ minimizes $f(\cdot,Z_k)$, one also has
\begin{align}
\label{eq:appendix_ss_grad_bound}
\|\nabla f(x_k,Z_k)\|^2 \leq 2L\,\Delta(x_k,Z_k).
\end{align}

\noindent We now consider the three step-size policies separately.

\noindent (i): Suppose $\eta_k=\eta\in(0,2/L)$ for every $k\geq 1$. 

\begin{itemize}
\item If $\eta\in(0,1/L]$, then the last term in eq.~\eqref{eq:appendix_ss_master} is non-positive, and therefore
\begin{align*}
\|x_{k+1}-x_\star\|^2
\le
\|x_k-x_\star\|^2-2\eta\,\Delta(x_k,Z_k).
\end{align*}
Hence Assumption~\ref{ass:ss_condition} holds with $\phi_k=\eta$ and $\bar\phi=\eta$.

\item If instead $\eta\in(1/L,2/L)$, then combining eq.~\eqref{eq:appendix_ss_master} with eq.~\eqref{eq:appendix_ss_grad_bound} yields
\begin{align*}
\|x_{k+1}-x_\star\|^2
&\le
\|x_k-x_\star\|^2
-2\eta\,\Delta(x_k,Z_k)
+2\eta L \Bigl(\eta-\frac{1}{L}\Bigr)\cdot \Delta(x_k,Z_k)\\
&=
\|x_k-x_\star\|^2
-2\eta(2-L\eta)\,\Delta(x_k,Z_k).
\end{align*}
Hence Assumption~\ref{ass:ss_condition} holds with $\phi_k=\eta(2-L\eta)$ and $\bar\phi=\eta(2-L\eta)$.
\end{itemize}
    
\noindent (ii): Suppose the step size $\eta_k$ is selected by the stochastic line-search rule in \Cref{alg:backtrack}. We start recalling that, by convexity of $f(\cdot, Z_k)$ (\Cref{ass:convexity}), it holds that
\begin{align}\label{eq:prop_examples_1}
\Delta(x_k, Z_k) = f(x_{k}, Z_k) - f(x_{\star}, Z_k) \leq \langle \nabla f(x_k, Z_k), x_k - x_{\star} \rangle.
\end{align}
Second, the line search condition in line 4 of \Cref{alg:backtrack} ensures that 
\begin{align*}
f(x_{k+1},Z_k)\le f(x_k,Z_k)-c\eta_k\|\nabla f(x_k,Z_k)\|^2,    
\end{align*}
which, by interpolation (\Cref{ass:interpolation}), implies
\begin{align}\label{eq:prop_examples_2}
c\eta_k\|\nabla f(x_k,Z_k)\|^2\le \Delta(x_k,Z_k).
\end{align}
Substituting eq.~\eqref{eq:prop_examples_1} and eq.~\eqref{eq:prop_examples_2} into eq.~\eqref{eq:appendix_ss_general_expansion} gives
\begin{align*}
\|x_{k+1}-x_\star\|^2
&\le
\|x_k-x_\star\|^2
-2\eta_k\,\Delta(x_k,Z_k)
+\eta_k^2\|\nabla f(x_k,Z_k)\|^2\\
&\le
\|x_k-x_\star\|^2
-2\eta_k\,\Delta(x_k,Z_k)
+\frac{\eta_k}{c}\,\Delta(x_k,Z_k)\\
&=
\|x_k-x_\star\|^2
-\frac{2c-1}{c}\eta_k\,\Delta(x_k,Z_k).
\end{align*}
\noindent Hence Assumption~\ref{ass:ss_condition} holds with
\begin{align*}
\phi_k=\frac{2c-1}{2c}\eta_k, \quad \bar\phi=\frac{2c-1}{2c}\bar\eta,
\end{align*}
where the choice of $\bar{\phi}$ is due to \Cref{lemma:sls}.

\noindent (iii): Suppose
\begin{align*}
\eta_k=
\begin{cases}
\dfrac{\Delta(x_k,Z_k)}{\|\nabla f(x_k,Z_k)\|^2}, & \text{if } \nabla f(x_k,Z_k)\neq 0,\\[3mm]
0, & \text{otherwise}.
\end{cases}
\end{align*}
If $\nabla f(x_k,Z_k)=0$, then by convexity of $f(\cdot,Z_k)$ (eq.~\eqref{eq:prop_examples_1}) we have $\Delta(x_k,Z_k)=0$, and therefore
\[
\|x_{k+1}-x_\star\|^2=\|x_k-x_\star\|^2
\le
\|x_k-x_\star\|^2-2\cdot \frac{1}{4L}\cdot \Delta(x_k,Z_k).
\]
If instead $\nabla f(x_k,Z_k)\neq 0$, then using convexity of $f(\cdot, Z_k)$ in eq.~\eqref{eq:appendix_ss_general_expansion},
\begin{align*}
\|x_{k+1}-x_\star\|^2
&\le
\|x_k-x_\star\|^2
-2\eta_k\Delta(x_k,Z_k)
+\eta_k^2\|\nabla f(x_k,Z_k)\|^2\\
&=
\|x_k-x_\star\|^2
-\frac{\Delta(x_k,Z_k)^2}{\|\nabla f(x_k,Z_k)\|^2}.
\end{align*}
Using eq.~\eqref{eq:appendix_ss_grad_bound}, we obtain
\[
\frac{\Delta(x_k,Z_k)^2}{\|\nabla f(x_k,Z_k)\|^2}
\ge
\frac{1}{2L}\Delta(x_k,Z_k),
\]
and therefore
\[
\|x_{k+1}-x_\star\|^2
\le
\|x_k-x_\star\|^2-\frac{1}{2L}\Delta(x_k,Z_k)
=
\|x_k-x_\star\|^2-2\cdot \frac{1}{4L}\Delta(x_k,Z_k).
\]
Hence Assumption~\ref{ass:ss_condition} holds with
\[
\phi_k=\frac{1}{4L},
\qquad
\bar\phi=\frac{1}{4L}.
\]
This concludes the proof.

\subsection{Missing Proofs for Theorem \ref{thm:main_result}}
\subsubsection{Proof of Lemma \ref{lemma:mds_prop}}\label{subsubsec:appendix_convex_smooth_optimization_proof_lemma_mds_prop}
\ref{lemma:mds_boundedness}. Note that, for every $k \in \bbN$, it holds that
\begin{align}
\label{eq:lemma_mds_prop_boundedness_1}
0 \leq \Delta(x_k, Z_k) \overset{(a)}{\leq} \frac{L \cdot \|x_k - x_{\star}\|^2}{2} \overset{(b)}{\leq} \frac{L \cdot \|x_1 - x_{\star}\|^2}{2} \quad \text{(pathwise)},
\end{align}
where: (a) follows from $f(\cdot, Z_k)$ being $L$-smooth (\Cref{ass:convexity}) and the optimality of $x_{\star}$ (\Cref{ass:interpolation}), and (b) follows from  eq.~\eqref{eq:fejer_monotonicity}. Similarly, it holds that
\begin{align*}
0 &\leq \Delta(x_k) \overset{(a)}{\leq} \frac{L \cdot \|x_k - x_{\star}\|^2}{2} \overset{(b)}{\leq} \frac{L \cdot \|x_1 - x_{\star}\|^2}{2} \quad \text{(pathwise)},
\end{align*}
where: (a) follows from $F$ being $L$-smooth (\Cref{ass:convexity}), and (b) follows from eq.~\eqref{eq:fejer_monotonicity}.
Therefore, it holds that
\begin{align*}
|\zeta_k| &= |\Delta(x_k) - \Delta(x_k, Z_k)| \leq \frac L 2 \|x_1 - x_{\star}\|^2.
\end{align*}    
\noindent \ref{lemma:mds}. First note that, by the above item, it holds that
\begin{align*}
\mathbb{E}[|\zeta_k|] &\leq \frac L 2 \|x_1 - x_{\star}\|^2 < +\infty.
\end{align*}
Moreover,
\begin{align*}
\mathbb{E}_k[\zeta_k] &= \mathbb{E}_k[\Delta(x_k) - \Delta(x_k, Z_k)] \nonumber\\
&= \mathbb{E}_k[\Delta(x_k)] - \mathbb{E}_k[\Delta(x_k, Z_k)] \nonumber\\
&\overset{(a)}{=} \Delta(x_k) - \mathbb{E}_k[\Delta(x_k, Z_k)] \nonumber\\
&= (F(x_k) - F(x_{\star})) - \mathbb{E}_k[f(x_k, Z_k) - f(x_{\star}, Z_k)] \nonumber\\
&= F(x_k) - F(x_{\star}) - \mathbb{E}_k[f(x_k, Z_k)] + \mathbb{E}_k[f(x_{\star}, Z_k)] \nonumber\\
&= F(x_k) - F(x_{\star}) - F(x_k) + F(x_{\star}) \nonumber\\
&= 0,
\end{align*}
where (a) is due to $x_k$ being $\sigma(Z_1,\dots,Z_{k-1})$-measurable. Therefore $(\zeta_k)_{k \geq 1}$ is a martingale difference sequence.

\noindent \ref{lemma:mds_conditional_variance}. The item follows by the following elementary calculations
\begin{align*}
\mathbb{E}_k[\zeta_k^2] &= \mathbb{E}_k[(\Delta(x_k) - \Delta(x_k, Z_k))^2] \\
&= \mathbb{E}_k[\Delta(x_k, Z_k)^2 -2\Delta(x_k, Z_k)\Delta(x_k) + \Delta(x_k)^2] \\
&= \mathbb{E}_k[\Delta(x_k, Z_k)^2 -2\mathbb{E}_k[\Delta(x_k, Z_k)\Delta(x_k)] + \mathbb{E}_k[\Delta(x_k)^2] \\
&\overset{(a)}{=} \mathbb{E}_k[\Delta(x_k, Z_k)^2] -2\mathbb{E}_k[\Delta(x_k, Z_k)] \Delta(x_k) + \Delta(x_k)^2 \\
&= \mathbb{E}_k[\Delta(x_k, Z_k)^2 -2\Delta(x_k)^2 + \Delta(x_k)^2 \\
&= \mathbb{E}_k[\Delta(x_k, Z_k)^2] - \Delta(x_k)^2 \\
&\leq  \mathbb{E}_k[\Delta(x_k, Z_k)^2] \\
&= \mathbb{E}_k[\Delta(x_k, Z_k) \cdot \Delta(x_k, Z_k)] \\
&\overset{(b)}{\leq} \frac L 2 \|x_1 - x_{\star}\|^2 \cdot \mathbb{E}_k[\Delta(x_k, Z_k)] \\
&= \frac L 2 \|x_1 - x_{\star}\|^2 \cdot \Delta(x_k) \quad \text{(pathwise),}
\end{align*}
where: (a) is due to $x_k$ being $\sigma(Z_1,\dots,Z_{k-1})$-measurable, and (b) follows from eq.~\eqref{eq:lemma_mds_prop_boundedness_1}. This concludes the proof.

\subsubsection{Proof of Lemma \ref{lemma:bernstein_condition}}\label{subsubsec:appendix_convex_smooth_optimization_proof_lemma_bernstein_condition}
We start noting that for all $\lambda \geq 0$, by virtue of eq.~\eqref{eq:mds_boundedness}, it holds $\lambda \cdot \zeta_k \leq  \lambda \cdot \nu$ pathwise. Moreover, by elementary analysis, one has
\begin{align*}
\forall x < 3 \quad \exp(x) \leq 1 + x + \frac{x^2}{2 \cdot (1-x/3)}.
\end{align*}
Combining the above observations yields
\begin{align*}
\forall \lambda \in [0,3/\nu) \quad \exp(\lambda \cdot \zeta_k) &\leq 1 + \lambda \cdot \zeta_k + \frac{\lambda^2 \cdot \zeta_k^2}{2 \cdot (1-\lambda \cdot \zeta_k/3)} \\
&\leq 1 + \lambda \cdot \zeta_k + \frac{\lambda^2 \cdot \zeta_k^2}{2 \cdot (1-\lambda \cdot \nu/3)}.
\end{align*}
Taking conditional expectation on both sides of the above display, one obtains that (pathwise) the following holds
\begin{align*}
\forall \lambda \in [0,3/\nu) \quad \bbE_k[\exp(\lambda \cdot \zeta_k)] &\leq \bbE_k\bigg[1 + \lambda \cdot \zeta_k + \frac{\lambda^2 \cdot \zeta_k^2}{2 \cdot (1-\lambda \cdot \nu/3)}\bigg] \\
&\overset{(a)}{=} 1 + \lambda \cdot \bbE_k[\zeta_k] + \lambda^2 \cdot \frac{\bbE_k[\zeta_k^2]}{2\cdot(1-\lambda \cdot \nu/3)} \\
&\overset{(b)}{=} 1 + \lambda^2 \cdot \frac{\bbE_k[\zeta_k^2]}{2 \cdot (1-\lambda \cdot \nu/3)}  \\
&\overset{(c)}{\leq} 1 + \frac{\lambda^2 \cdot \nu}{2 \cdot (1-\lambda \cdot \nu/3)} \cdot \Delta(x_k)  \\
&\overset{(d)}{\leq} \exp \bigg(\frac{\lambda^2 \cdot \nu}{2 \cdot (1-\lambda \cdot \nu/3)} \cdot \Delta(x_k) \bigg),
\end{align*}
where: (a) follows from linearity of (conditional) expectation, (b) follows from \Cref{lemma:mds_prop}~\cref{lemma:mds}, (c) follows from \Cref{lemma:mds_prop}~\cref{lemma:mds_conditional_variance}, and (d) follows from $1+x \leq \exp(x)$ for every $x \in \mathbb{R}$. This concludes the proof.

\newpage
\section{Appendix: Section \ref{sec:sc}}\label{sec:appendix_sc}

\subsection{Full Statement of \Cref{thm:markov_p}}\label{appendix:sc_ub_full}
The result of this section is the following.
\begin{theorem}
\label{thm:markov_p_full}
Suppose Assumptions~\ref{ass:feasibility}, \ref{ass:interpolation}, and \ref{ass:strong_convexity} hold. Furthermore, suppose \Cref{ass:ss_condition} with $\bar{\phi} \leq 1/L$. Let $x_1 \in \mathcal{X}$ and $\kappa = L/\mu$, then the iterates generated by eq.~\eqref{eq:sgd_ls} satisfy:
\begin{enumerate}[label={\rm(\roman*)}]
\item For every $n \geq 1$, it holds that
\begin{align}
\label{eq:sc_in_expectation_full}
\bbE[\Delta(x_{n+1})] &\leq \frac{L \cdot \|x_{1} -x_{\star}\|^2}{2} \cdot \bigg(1- \bar{\phi} \cdot L \cdot \kappa^{-1} \bigg)^n.
\end{align}
\item For every $n \geq 1$, $\beta \in (0,1)$, and $p \in (0,1]$ it holds that
\begin{align}
\label{eq:hp_markov_p_full}
\bbP \bigg(\Delta(x_{n+1}) \leq \frac{1}{\beta^p} \cdot \frac{L \cdot \|x_{1} -x_{\star}\|^2}{2} \cdot \bigg(1- a_p(\bar{\phi}, L) \cdot \bar{\phi} \cdot L \cdot \kappa^{-1} \bigg)^{pn} \bigg) \geq 1-\beta.
\end{align}
where
\begin{align*}
a_p(\bar{\phi}, L) \coloneqq \frac{1-\big(1 - L \cdot \bar{\phi}\big)^{1/p}}{L \cdot \bar{\phi}}.
\end{align*}
\end{enumerate}
\end{theorem}
We note that \Cref{eq:sc_in_expectation_full} is not the main novelty of the present section, as similar in-expectation guarantees have appeared in prior works, although only in the more specific setting of finite-sum optimization and for specific step-size policies. We report it here primarily to put the high-probability result into perspective, and for its slightly greater generality.

\subsection{Proof of Theorem \ref{thm:markov_p_full}}\label{subsec:appendix_sc_thm_markov_p}
\rm{(i)}: Recalling eq.~\eqref{eq:one_step_inequality}, for every $k \geq 1$, it holds that
\begin{align}\label{eq:thm_sc_markov_p_eq1}
\|x_{k+1}-x_{\star}\|^2 &\leq \|x_{k} -x_{\star}\|^2 - 2 \cdot \phi_k \cdot \Delta(x_k, Z_k) \nonumber\\
&\leq \|x_{k} -x_{\star}\|^2 - 2 \cdot \bar{\phi} \cdot \Delta(x_k, Z_k),
\end{align}
where the second inequality follows from \Cref{ass:ss_condition} and $\Delta(x_k, Z_k) \geq 0$. Thus, taking conditional expectation on both sides of the above display yields 
\begin{align}
\label{eq:conditional_expctation_sc}
\bbE_{k}[\|x_{k+1}-x_{\star}\|^2] &\leq \bbE_{k}\bigg[\|x_{k} -x_{\star}\|^2 - 2 \cdot \bar{\phi} \cdot \Delta(x_k, Z_k)\bigg] \nonumber\\
&\overset{(a)}{=} \|x_{k} -x_{\star}\|^2 - 2 \cdot \bar{\phi} \cdot \bbE_{k}\big[\Delta(x_k, Z_k)\big] \nonumber \\
&\overset{(b)}{=} \|x_{k} -x_{\star}\|^2 - 2 \cdot \bar{\phi} \cdot \Delta(x_k) \nonumber \\
&\overset{(c)}{\leq} \|x_{k} -x_{\star}\|^2 - \bar{\phi} \cdot \mu \cdot \|x_k - x_{\star}\|^2 \nonumber \\
&= \bigg(1 - \bar{\phi} \cdot \mu  \bigg) \cdot \|x_{k} -x_{\star}\|^2. 
\end{align}
where: (a) follows from the linearity of (conditional) expectation and $\|x_k-x_{\star}\|^2$ being $\sigma(Z_1,\dots,Z_{k-1})$-measurable, (b) follows from eq.~\eqref{eq:cond_exp}, and (c) follows from the strong convexity of $F$ (\Cref{ass:strong_convexity}).

\noindent Using the tower rule of expectation, unrolling eq.~\eqref{eq:conditional_expctation_sc}, and noting that $\|x_1-x_{\star}\|^2$ is fixed, yields 
\begin{align}
\label{eq:iterates_convergence_sc_inexpectation}
\bbE[\|x_{n+1}-x_{\star}\|^2] &\leq \bigg(1- \bar{\phi} \cdot \mu \bigg)^n \cdot \|x_{1} -x_{\star}\|^2.
\end{align}
Finally, combining the above with the $L$-smoothness of $F$ (\Cref{ass:convexity}), yields
\begin{align}
\label{eq:sc_in_expectation_results_1}
\bbE[\Delta(x_{n+1})] &\leq \frac{L \cdot \|x_{1} -x_{\star}\|^2}{2} \cdot \bigg(1- \bar{\phi} \cdot \mu \bigg)^n.
\end{align}
\noindent This proves the first claim. 

\noindent \rm{(ii)}: For $p\in(0,1]$, let us define
\begin{align*}
a_p(\bar{\phi}, L) \coloneqq \frac{1-\big(1-L \cdot \bar{\phi}\big)^{1/p}}{L \cdot \bar{\phi}},
\end{align*}
and note that $a_p(\bar \phi, \eta) \geq0$ since $L \cdot \bar{\phi} \leq 1$ by assumption. Recall that eq.~\eqref{eq:thm_sc_markov_p_eq1}, for every $k \geq 1$, yields 
\begin{align*}
\|x_{k+1}-x_{\star}\|^2 \leq \|x_{k} -x_{\star}\|^2 - 2 \cdot \bar{\phi} \cdot \Delta(x_k, Z_k).
\end{align*}

\noindent For every $k \geq 1$, define
\begin{align*}
u_k
\coloneqq
\begin{cases}
\displaystyle \frac{2 \cdot \bar{\phi} \cdot \Delta(x_k,Z_k)}{\|x_k-x_\star\|^2},
& \text{if } x_k \neq x_\star,\\[3mm]
0, & \text{if } x_k=x_\star.
\end{cases}
\end{align*}
Then
\begin{align*}
\|x_{k+1}-x_\star\|^2
\leq
\|x_k-x_\star\|^2 \cdot (1-u_k).
\end{align*}
Moreover, since $f(\cdot, Z_k)$ is $L$-smooth (\Cref{ass:strong_convexity}) and minimized at $x_\star$ (\Cref{ass:interpolation}), we have
\begin{align*}
\Delta(x_k,Z_k)\leq \frac{L \cdot \|x_k-x_\star\|^2}{2},
\end{align*}
and therefore
\begin{align*}
0\leq u_k\leq L \cdot \bar{\phi}.
\end{align*}
Now consider the function
\begin{align*}
\varphi(u)\coloneqq (1-u)^{1/p}
\end{align*}
on the interval $[0, L \cdot \bar{\phi}]$. Since $1/p\geq 1$, the function $\varphi$ is convex. Hence, for every $u$ in this interval, a convex function is below the secant, that is
\begin{align}\label{eq:thm_sc_markov_p_eq2}
(1-u)^{1/p} &\leq 1 - \frac{1-\big(1-L \cdot \bar{\phi}\big)^{1/p}}{L \cdot \bar{\phi}} \cdot u \nonumber \\ 
&= 1 - a_p(\bar{\phi}, L) \cdot u.
\end{align}
Applying this with $u=u_k$, we obtain
\begin{align*}
\|x_{k+1}-x_\star\|^{2/p} &\leq \|x_k-x_\star\|^{2/p} \cdot (1-u_k)^{1/p} \\
&\overset{(a)}{\leq} \|x_k-x_\star\|^{2/p} - a_p(\bar{\phi}, L) \cdot u_k \cdot \|x_k-x_\star\|^{2/p} \\
&= \|x_k-x_\star\|^{2/p} - 2 \cdot a_p(\bar{\phi}, L) \cdot \bar{\phi} \cdot \|x_k-x_\star\|^{2/p-2} \cdot \Delta(x_k,Z_k),
\end{align*}
where $(a)$ follows from eq.~\eqref{eq:thm_sc_markov_p_eq2}. Taking conditional expectation on both sides of the above display yields
\begin{align*}
\bbE_k [\|x_{k+1}-x_\star\|^{2/p}] &\overset{(a)}{\leq} \|x_k-x_\star\|^{2/p} - 2 \cdot a_p(\bar{\phi}, L) \cdot \bar{\phi} \cdot \|x_k-x_\star\|^{2/p-2} \cdot \bbE_k [ \Delta(x_k,Z_k)] \\ 
&= \|x_k-x_\star\|^{2/p} - 2 \cdot a_p(\bar{\phi}, L) \cdot \bar{\phi} \cdot \|x_k-x_\star\|^{2/p-2} \cdot \Delta(x_k) \\
&\overset{(b)}{\leq} \|x_k-x_\star\|^{2/p} - a_p(\bar{\phi}, L) \cdot \bar{\phi} \cdot \mu \cdot \|x_k-x_\star\|^{2/p}\\
&= (1 - a_p(\bar{\phi}, L) \cdot \bar{\phi} \cdot \mu) \cdot \|x_k-x_\star\|^{2/p},
\end{align*}
where (a) follows from the linearity of (conditional) expectation and $\|x_k-x_{\star}\|^{2/p}$ being $\sigma(Z_1,\dots,Z_{k-1})$-measurable, and (b) follows by the strong convexity of $F$ (\Cref{ass:strong_convexity}). Using the tower rule of expectation and unrolling the above recursion, we get
\begin{align*}
\bbE[\|x_{n+1}-x_\star\|^{2/p}] \leq (1 - a_p(\bar{\phi}, L) \cdot \bar{\phi} \cdot \mu)^n \cdot  \|x_1-x_\star\|^{2/p}.
\end{align*}
Therefore,
\begin{align}
\label{eq:p_moment_bound}
\bbE[\Delta(x_{n+1})^{1/p}] &\overset{(a)}{\leq} \bigg(\frac{L}{2}\bigg)^{1/p}\bbE[\|x_{n+1}-x_\star\|^{2/p}] \nonumber \\
&\leq \bigg(\frac{L\cdot \|x_1-x_\star\|^2}{2}\bigg)^{1/p} \cdot (1 - a_p(\bar{\phi}, L) \cdot \bar{\phi} \cdot \mu)^n,
\end{align}
where (a) follows from the $L$-smoothness of $F$ (\Cref{ass:strong_convexity}). 

\noindent Now define
\begin{align*}
t &\coloneqq \frac{1}{\beta^p} \cdot \frac{L\cdot \|x_1-x_\star\|^2}{2} \cdot (1-a_p(\bar{\phi},L) \cdot \bar{\phi} \cdot \mu)^{pn} \\
&= \frac{1}{\beta^p} \cdot \frac{L\cdot \|x_1-x_\star\|^2}{2} \cdot (1-a_p(\bar{\phi},L) \cdot \bar{\phi} \cdot L \cdot \kappa^{-1})^{pn}.
\end{align*}
Therefore, by Markov's inequality and eq.~\eqref{eq:p_moment_bound},
\begin{align*}
\bbP\big(\Delta(x_{n+1})\ge t\big)
&=
\bbP\big(\Delta(x_{n+1})^{1/p}\ge t^{1/p}\big)\\
&\le
\frac{\bbE[\Delta(x_{n+1})^{1/p}]}{t^{1/p}}
\le
\beta.
\end{align*}
Equivalently,
\begin{align*}
\bbP\big(\Delta(x_{n+1})\le t\big)\ge 1-\beta,
\end{align*}
which proves the second claim.

\subsection{Proof of Theorem \ref{thm:sc_lower_bound_p}}\label{subsec:appendix_sc_lower_bound_p}
Fix any $p\in(0,1]$ and any $\kappa>1$. Consider the one-dimensional problem defined by
\[
\mathcal{X}=\mathbb{R},
\qquad
f(x,Z)\coloneqq \frac{Z}{2} \cdot x^2,
\]
where the random variable $Z$ takes values in $\{0,L\}$ and satisfies
\[
\mathbb{P}(Z=0)=1-\frac{1}{\kappa},
\qquad
\mathbb{P}(Z=L)=\frac{1}{\kappa}.
\]
Then, for every realization of $Z$, the function $f(\cdot,Z)$ is convex and $L$-smooth. Moreover, $x_\star=0$ minimizes $f(\cdot,Z)$ for every $Z\in\{0,L\}$, so Assumption~\ref{ass:interpolation} holds. Finally,
\[
F(x)=\mathbb{E}[f(x,Z)]
=
\frac{1}{2}\cdot \frac{L}{\kappa} \cdot x^2,
\]
hence $F$ is strongly convex with modulus $\mu = L/\kappa$ and condition number $\kappa$.

\noindent Now for $n \geq 1$, consider the event  
\[
E_n = \{Z_1=\cdots=Z_n=0\},
\]
and note that, by the independence of the $(Z_k)_{k \in [n]}$, it holds that
\begin{equation}
\label{eq:lb_event_probability}
\mathbb{P}(E_n)=\bigg(1-\frac{1}{\kappa}\bigg)^n.
\end{equation}
Let $(x_k)_{k\ge 1}$ be the iterates generated by eq.~\eqref{eq:sgd_ls} starting from $x_1\in\mathcal{X}$ with $\Delta(x_1)>0$. Since $\nabla f(x,Z) = Z \cdot x$, on $E_n$ we have for every $1\leq k\leq n$, $\nabla f(x_k,Z_k)=0$.
Thus 
\begin{equation}
\label{eq:lb_no_progress}
x_{n+1}=x_1
\qquad\text{and}\qquad
\Delta(x_{n+1})=\Delta(x_1).
\end{equation}
\noindent Fix now any
\[
a>-p\log\Bigl(1-\frac1\kappa\Bigr).
\]
Thus, for every $n\ge 1$,
\[
e^{-an/p}<\bigg(1-\frac1\kappa\bigg)^n,
\]
so that we can pick $\beta\in(0,1)$ s.t.
\begin{equation}
\label{eq:lb_choice_beta}
e^{-an/p}<\beta<\bigg(1-\frac1\kappa\bigg)^n.
\end{equation}

\noindent Now, by the left-hand side of eq.~\eqref{eq:lb_choice_beta},
\[
\beta^{-p}e^{-an}<1.
\]
Therefore, on the event $E_n$, using eq.~\eqref{eq:lb_no_progress}, we have
\[
\Delta(x_{n+1})
=
\Delta(x_1)
>
\beta^{-p}e^{-an}\Delta(x_1).
\]
Thus the event
\[
\bigg\{
\Delta(x_{n+1})
\le
\beta^{-p}e^{-an}\Delta(x_1)
\bigg\}
\]
is disjoint from $E_n$. Consequently,
\[
\mathbb{P}\bigg(
\Delta(x_{n+1})
\le
\beta^{-p}e^{-an}\Delta(x_1)
\bigg)
\le
1-\mathbb{P}(E_n).
\]
Using eq.~\eqref{eq:lb_event_probability} and the right-hand side of eq.~\eqref{eq:lb_choice_beta}, we conclude that
\[
\mathbb{P}\bigg(
\Delta(x_{n+1})
\le
\beta^{-p}e^{-an}\Delta(x_1)
\bigg)
\le
1-\bigg(1-\frac1\kappa\bigg)^n
<
1-\beta.
\]
This proves the claim.

\section{Appendix: Section \ref{sec:extensions}}\label{appendix:extensions}
Throughout this section we consider SGD with constant step size as given by the following. For $x_1 \in \cX$ and $\eta \geq 0$, we consider the following recurrence
\begin{align}
\label{eq:sgd_cs}
\forall k \geq 1, \quad x_{k+1} = x_k - \eta \cdot \nabla f(x_k,Z_k).
\end{align}

\subsection{Low Noise Regime}\label{appendix:extensions_lownoise}
The result of this section is the following.
\begin{theorem}\label{thm:hp_generalized_interpolation_full}
Assume Assumptions \ref{ass:feasibility}, \ref{ass:convexity}, and \ref{ass:exp_noise_optimum}, let $x_1 \in \cX$ and $\eta \in (0, 2/L)$. Then the iterates generated by eq. \eqref{eq:sgd_cs} satisfy the following. For every $n \geq 1$ and every $\beta \in (0,1)$, with probability at least $1-\beta$, it holds that
\begin{align*}
\Delta(\bar x_n) &\leq \frac{4}{(2 - L \cdot \eta) \cdot \eta} \cdot \frac{\|x_1-x_\star\|^2}{n} + \frac{4}{3} \cdot \log\Big(\frac 3 \beta \Big) \cdot \frac{L \cdot  R_\eta^2}{n} \\
&\qquad + 16 \cdot \sqrt{2\log(2)} \cdot \frac{\sigma_\star \cdot R_\eta}{2 - L\cdot \eta} \cdot  \sqrt{\frac{\log\Big(\frac 3 \beta \Big)}{n}} + 8 \cdot \frac{\eta \cdot (2 + L \cdot \eta)}{(2 - L \cdot \eta)^2} \cdot  \sigma_\star^2 \cdot \log\Big(\frac 3 \beta \Big),
\end{align*}
where
\begin{align*}
R_\eta \coloneqq 2 \cdot \|x_1-x_\star\| + 6 \cdot \eta \cdot \sigma_\star \sqrt{\frac{2+L \cdot \eta}{2-L \cdot \eta} \cdot n \cdot \log\Big(\frac 3 \beta \Big)}.
\end{align*}
\noindent Moreover, if $x_1  = x_{\star}$ set $\eta = 0$, and then $\Delta(\bar x_{n}) = 0$ (pathwise). Otherwise, let
\begin{align*}
\eta = \min\left\{ \frac{1}{4 \cdot L}, \frac{\|x_1-x_\star\|}{32 \cdot \sigma_\star \cdot \sqrt{n \cdot \log(3/\beta)}} \right\}.
\end{align*}
Then, with probability at least $1-\beta$, it holds that
\[
\Delta(\bar x_n) 
\le
8 \cdot \bigg(2 + \log \Big(\frac 3 \beta\Big)\bigg) \cdot \frac{L \cdot \|x_1-x_\star\|^2}{n}
+
98 \cdot 
\sigma_\star \cdot \|x_1-x_\star\| \cdot 
\sqrt{\frac{\log\Big( \frac 3 \beta \Big)}{n}}.
\]
\end{theorem}

\subsubsection{Proof of Theorem \ref{thm:hp_generalized_interpolation_full}}
\paragraph{Additional notation.}
To control the location of the iterates, for every $R>0$ we define 
\[
\tau_R \coloneqq \inf\{k\ge 1:\|x_k-x_\star\|>R\}.
\]
The following two quantities play a key role in the analysis. For every $n\ge 1$, let
\begin{align*}
M_n &\coloneqq \sum_{k=1}^n \left\langle \nabla f(x_\star,Z_k), x_\star - x_k \right\rangle, \\
G_n &\coloneqq \sum_{k=1}^n \Big( \Delta(x_k) - \Delta(x_k,Z_k) + \langle \nabla f(x_\star,Z_k), x_k - x_\star \rangle \Big).
\end{align*}
For every $R>0$, define their truncated counterparts by
\begin{align*}
M_n^{(R)} &\coloneqq \sum_{k=1}^n \left\langle \nabla f(x_\star,Z_k), x_\star - x_k \right\rangle \mathbf 1_{\{k<\tau_R\}}, \\
G_n^{(R)} &\coloneqq \sum_{k=1}^n \Big( \Delta(x_k) - \Delta(x_k,Z_k) + \langle \nabla f(x_\star,Z_k), x_k - x_\star \rangle \Big)\mathbf 1_{\{k<\tau_R\}}.
\end{align*}

\paragraph{Preliminaries.} 
We start with the following lemmas, for which we assume always that Assumptions \ref{ass:feasibility}, \ref{ass:convexity}, and \ref{ass:exp_noise_optimum} hold. 
\begin{lemma}
\label{lem:localized_Mn_bound}
For every $R>0$, every $n\ge 1$, and every $\beta\in(0,1)$, it holds that
\begin{align}
\bbP\left(\sup_{m\in[n]} M_m^{(R)} \geq 2 \cdot \sigma_\star \cdot R \cdot \sqrt{2 \cdot \log(2) \cdot n \cdot \log(3/\beta)} \right)
\leq \frac{\beta}{3}.
\end{align}
\end{lemma}

\begin{proof}
For every $k\in[n]$, define
\[
X_k \coloneqq \left\langle \nabla f(x_\star,Z_k), x_\star - x_k \right\rangle \mathbf 1_{\{k<\tau_R\}}.
\]
Since $x_k$ and $\mathbf 1_{\{k<\tau_R\}}$ are $\sigma(Z_1, \dots, Z_{k-1})$-measurable, and $Z_k$ is independent from the past, $\bbE_k[\nabla f(x_{\star}, Z_k)] = \nabla F(x_{\star}) = 0$, and by \Cref{ass:exp_noise_optimum}, \Cref{lem:exp_square_implies_linear_subgaussian} applied conditionally with
\[
X=\nabla f(x_\star,Z_k) \qquad\text{and}\qquad u=(x_\star - x_k)\mathbf 1_{\{k<\tau_R\}}
\]
implies that, for every $\lambda\in\bbR$, almost surely,
\begin{align*}
\bbE_k\left[ \exp(\lambda X_k) \right] &\leq \exp\left( 2\log(2)\,\lambda^2\sigma_\star^2 \|(x_\star-x_k)\mathbf 1_{\{k<\tau_R\}}\|^2 \right) \\
&= \exp\left( 2\log(2)\,\lambda^2\sigma_\star^2  \|x_\star-x_k\|^2\mathbf 1_{\{k<\tau_R\}} \right) \\
&\leq \exp\left( 2\log(2)\,\lambda^2\sigma_\star^2  R^2 \right).    
\end{align*}
For any $\lambda>0$, define
\[
\Psi_0(\lambda)\coloneqq 1, \qquad \Psi_m(\lambda) \coloneqq \exp\left(
\lambda \sum_{k=1}^m X_k - 2\log(2)\,\lambda^2\sigma_\star^2R^2 m \right), \qquad \forall m\in[n].
\]
Then $(\Psi_m(\lambda))_{m=0}^n$ is a nonnegative supermartingale. Indeed,
\[
\bbE_m[\Psi_{m+1}(\lambda)] = \Psi_m(\lambda) \bbE_m\left[ \exp\left( \lambda X_{m+1} - 2\log(2) \lambda^2\sigma_\star^2R^2 \right) \right] \leq \Psi_m(\lambda).
\]
Hence, by \Cref{lem:doob_maximal}, for every $v>0$,
\begin{align}\label{eq:doob_application}
\bbP\left( \max_{0 \leq m \leq n} \Psi_m(\lambda)\ge v \right) \leq \frac{\bbE[\Psi_0(\lambda)]}{v} \leq \frac{1}{v}.
\end{align}
Note that for $t > 0$, if we let $v = \exp(\lambda t - 2\log(2)\,\lambda^2\sigma_\star^2R^2 n)$, then
\[
\max_{0 \leq m \leq n} \sum_{k=1}^m X_k \ge t \Rightarrow \max_{0 \leq m \leq n} \Psi_m(\lambda)\ge v,
\]
therefore, from eq. \eqref{eq:doob_application}, we get
\[
\bbP\left( \max_{m\in[n]} \sum_{k=1}^m X_k \geq t \right) \leq \exp\left( -\lambda t + 2\log(2)\,\lambda^2\sigma_\star^2R^2 n \right).
\]
Minimizing the right-hand side over $\lambda>0$ gives
\[
\bbP\left( \max_{m\in[n]} \sum_{k=1}^m X_k \geq t \right) \leq \exp\left( -\frac{t^2}{8\log(2)\,\sigma_\star^2R^2 n} \right).
\]
Choosing
\[
t = \sigma_\star R\sqrt{8\log(2)\,n\log(3/\beta)}
\]
concludes the proof.
\end{proof}
\begin{lemma}
\label{lem:quadratic_noise_bound}
For every $n \geq 1$ and every $\beta\in(0,1)$, it holds that
\begin{align*}
\bbP\left( \sum_{k=1}^n \|\nabla f(x_\star,Z_k)\|^2 \geq \sigma_\star^2 \cdot (\log(2) \cdot n + \log(3/\beta)) \right) \leq \frac{\beta}{3}.    
\end{align*}
\end{lemma}

\begin{proof}
By \Cref{ass:exp_noise_optimum} and the independence,
\[
\bbE\!\left[
\exp\!\left(
\frac{1}{\sigma_\star^2}
\sum_{k=1}^n
\|\nabla f(x_\star,Z_k)\|^2
\right)
\right]
=
\prod_{k=1}^n
\bbE\!\left[
\exp\!\left(
\frac{\|\nabla f(x_\star,Z_k)\|^2}{\sigma_\star^2}
\right)
\right]
\le
2^n.
\]
Hence, by Markov's inequality, for every $t>0$,
\[
\bbP\!\left(
\sum_{k=1}^n \|\nabla f(x_\star,Z_k)\|^2 \ge t
\right)
\le
2^n e^{-t/\sigma_\star^2}.
\]
Choosing
\[
t=\sigma_\star^2\bigl(n\log(2)+\log(3/\beta)\bigr)
\]
gives
\[
\bbP\!\left(
\sum_{k=1}^n \|\nabla f(x_\star,Z_k)\|^2 \ge t
\right)
\le
2^n e^{-n\log(2)}e^{-\log(3/\beta)}
=
\frac{\beta}{3}.
\]
This concludes the proof.
\end{proof}

\begin{lemma}
\label{lem:localized_Gn_bound}
For every $R > 0$, $ n \geq 1$ and every $\beta\in(0,1)$, it holds that 
\begin{align*}
\bbP \bigg( G_n^{(R)} \geq 
\frac 1 2 \cdot 
\sum_{k=1}^n \Delta(x_k) \cdot \mathbf 1_{\{k<\tau_R\}} +
\frac{2}{3} \cdot L\cdot R^2 \cdot \log(3/\beta) \bigg) \leq \frac{\beta}{3}.
\end{align*}
\end{lemma}

\begin{proof}
For every $k\in[n]$, define
\[
Y_k \coloneqq \Big( \Delta(x_k) - \Delta(x_k,Z_k) + \langle \nabla f(x_\star,Z_k),x_k - x_\star\rangle \Big)\mathbf 1_{\{k<\tau_R\}}.
\]
\noindent We note that $(Y_k)_{k=1}^n$ is a martingale difference sequence. Indeed, $\mathbf 1_{\{k<\tau_R\}}$ and $x_k$ are $\sigma(Z_1,\dots, Z_{k-1})$-measurable, and since $\nabla F(x_\star)=0$, we have $\bbE_k[Y_k]=0$. Moreover, by convexity and smoothness of $f(\cdot,Z_k)$ (\Cref{ass:convexity}),
\[
0 \leq \Delta(x_k,Z_k)-\langle \nabla f(x_\star,Z_k),x_k-x_\star\rangle \leq \frac{L}{2}\|x_k-x_\star\|^2,
\]
and since $F$ is also convex and $L$-smooth (\Cref{ass:convexity}) and $x_\star\in X_\star$,
\[
0\leq \Delta(x_k)\leq \frac{L}{2}\|x_k-x_\star\|^2.
\]
Hence, on the event $\{k<\tau_R\}$,
\[
0 \leq \Delta(x_k,Z_k)-\langle \nabla f(x_\star,Z_k),x_k-x_\star\rangle \leq
\frac{L R^2}{2}
\]
and
\[
0\le \Delta(x_k)\le \frac{L R^2}{2}.
\]
Therefore, $|Y_k| \leq (L R^2)/2$.

\noindent We now control the conditional variance. 
\begin{align*}
\bbE_k[Y_k^2] &= \mathrm{Var}_k\left( \big(\Delta(x_k,Z_k)-\langle \nabla f(x_\star,Z_k),x_k-x_\star\rangle\big)\mathbf 1_{\{k<\tau_R\}} \right) \\
&\le \bbE_k\left[ \big( \Delta(x_k,Z_k)-\langle \nabla f(x_\star,Z_k),x_k-x_\star\rangle \big)^2 \mathbf 1_{\{k<\tau_R\}} \right] \\
&\leq \frac{LR^2 }{2} \bbE_k\left[ \big( \Delta(x_k,Z_k)-\langle \nabla f(x_\star,Z_k),x_k-x_\star\rangle \big)\mathbf 1_{\{k<\tau_R\}}\right] \\
&= \frac{LR^2 }{2} \Delta(x_k) \mathbf 1_{\{k<\tau_R\}}.
\end{align*}
\noindent Thus we have showed that $(Y_k)_{k \in [n]}$ satisfies the same properties as the sequence $(\zeta_k)_{k \in [n]}$ in the proof of \Cref{thm:main_result}, with $\nu = (LR^2)/2$ and $R^2 \Delta(x_k) \mathbf 1_{\{k<\tau_R\}}$ in place of $\|x_1-x_{\star}\|^2 \Delta(x_k)$ in the conditional variance. Therefore, the very same argument used in that proof shows that 
\[
\bbP \bigg( G_n^{(R)} \geq 
\frac 1 2 \cdot 
\sum_{k=1}^n \Delta(x_k) \cdot \mathbf 1_{\{k<\tau_R\}} +
\frac{2}{3} \cdot L\cdot R^2 \cdot \log(3/\beta) \bigg) \leq \frac{\beta}{3},
\]
which proves the claim.
\end{proof}
\begin{lemma}[Error Decomposition]
\label{lemma:pathwise-low-noise}
Let $x_\star\in X_\star$ and $\eta \in (0, 2/L)$. Then, for every $n \geq 1$, it holds that
\begin{align*}
\sum_{k=1}^n \Delta(x_k) &\leq \frac{2}{(2 - L \cdot \eta) \cdot \eta} \cdot \|x_1-x_\star\|^2 + \frac{4}{2-L \cdot \eta} \cdot M_n \\
&\qquad + G_n + \frac{2\cdot \eta \cdot (2 + L\cdot \eta)}{(2 - L \cdot \eta)^2} \cdot \sum_{k=1}^n \|\nabla f(x_\star,Z_k)\|^2.
\end{align*}
\end{lemma}

\begin{proof}
For every $k \geq 1$, the updated rule \eqref{eq:sgd_ls} gives
\begin{align*}
\|x_{k+1}-x_\star\|^2
=
\|x_k-x_\star\|^2
-2\eta\left\langle \nabla f(x_k,Z_k),x_k-x_\star\right\rangle
+\eta^2\|\nabla f(x_k,Z_k)\|^2.    
\end{align*}
Summing from $k=1$ to $n$ and rearranging gives
\begin{align}
\label{eq:lemma_decomposition_eq1}
\sum_{k=1}^n \left\langle \nabla f(x_k,Z_k),x_k-x_\star\right\rangle \leq \frac{\|x_1-x_\star\|^2}{2\eta} + \frac{\eta}{2} \sum_{k=1}^n \|\nabla f(x_k,Z_k)\|^2.    
\end{align}
Now note that
\begin{align}
\label{eq:lemma_decomposition_eq2}
\|\nabla f(x_k,Z_k)\|^2 &\overset{(a)}{\leq} (1+\lambda)\|\nabla f(x_k,Z_k)-\nabla f(x_\star,Z_k)\|^2 + \left(1+\frac{1}{\lambda}\right)\|\nabla f(x_\star,Z_k)\|^2 \nonumber\\
&\overset{(b)}{\leq} L(1+\lambda) \left(\left\langle \nabla f(x_k,Z_k),x_k-x_\star\right\rangle - \left\langle \nabla f(x_\star,Z_k),x_k-x_\star\right\rangle \right) \nonumber \\
&\quad + \left(1+\frac{1}{\lambda}\right)\|\nabla f(x_\star,Z_k)\|^2,
\end{align}
where (a) holds for any $\lambda > 0$ by Young inequality \cite{Young1912}, and (b) follows from the convexity and the $L$-smoothness of $f(\cdot, Z_k)$ (\Cref{ass:convexity}). Therefore, plugging eq.~\eqref{eq:lemma_decomposition_eq2} into eq.~\eqref{eq:lemma_decomposition_eq1} yields
\begin{align*}
\sum_{k=1}^n \left\langle \nabla f(x_k,Z_k),x_k-x_\star\right\rangle &\leq \frac{\|x_1-x_\star\|^2}{2\eta} + \frac{L\eta(1+\lambda)}{2} \sum_{k=1}^n \left\langle \nabla f(x_k,Z_k),x_k-x_\star\right\rangle \\
&+ \frac{L\eta(1+\lambda)}{2} \cdot M_n + \frac{\eta}{2}\left(1+\frac{1}{\lambda}\right) \sum_{k=1}^n \|\nabla f(x_\star,Z_k)\|^2,
\end{align*}
and rearranging,
\begin{align*}
\left( 1-\frac{L\eta(1+\lambda)}{2} \right) \cdot \sum_{k=1}^n \left\langle \nabla f(x_k,Z_k),x_k-x_\star\right\rangle &\leq \frac{\|x_1-x_\star\|^2}{2\eta} + \frac{L\eta(1+\lambda)}{2} \cdot M_n \\
&\qquad + \frac{\eta}{2}\left(1+\frac{1}{\lambda}\right) \sum_{k=1}^n \|\nabla f(x_\star,Z_k)\|^2.    
\end{align*}
\noindent Now choose $\lambda = (2-L\eta)/(2L\eta)$, then the above becomes
\begin{align*}
\frac{2-L\eta}{4} \sum_{k=1}^n \left\langle \nabla f(x_k,Z_k),x_k-x_\star\right\rangle &\leq \frac{\|x_1-x_\star\|^2}{2\eta} + \frac{2+L\eta}{4} \cdot M_n \\
&\qquad + \frac{\eta(2+L\eta)}{2(2-L\eta)} \sum_{k=1}^n \|\nabla f(x_\star,Z_k)\|^2.
\end{align*}
Multiplying by $4/(2-L\eta)$, since $\eta \in (0,2/L)$, we get
\begin{align*}
\sum_{k=1}^n
\left\langle \nabla f(x_k,Z_k),x_k-x_\star\right\rangle \leq \frac{2\|x_1-x_\star\|^2}{(2-L\eta)\eta} + \frac{2+L\eta}{2 - L\eta} \cdot M_n + \frac{2\eta(2+L\eta)}{(2-L\eta)^2} \sum_{k=1}^n \|\nabla f(x_\star,Z_k)\|^2.
\end{align*}
\noindent On the other hand, using the convexity of $f(\cdot,Z_k)$ (\Cref{ass:convexity}) on the l.h.s. of the above, and the definition of $G_n$, yields
\begin{align*}
\sum_{k=1}^n \Delta(x_k) &\leq \frac{2\|x_1-x_\star\|^2}{(2-L\eta)\eta} + \left( \frac{2+L\eta}{2-L\eta}+1 \right)M_n + G_n + \frac{2\eta(2+L\eta)}{(2-L\eta)^2} \sum_{k=1}^n\|\nabla f(x_\star,Z_k)\|^2 \\
&= \frac{2\|x_1-x_\star\|^2}{(2-L\eta)\eta} + \frac{4}{2-L\eta} \cdot M_n + G_n + \frac{2\eta(2+L\eta)}{(2-L\eta)^2} \sum_{k=1}^n\|\nabla f(x_\star,Z_k)\|^2.
\end{align*}
This proves the claim.
\end{proof}

\paragraph{Main Argument.} We divide the proof into four steps.

\medskip
\noindent\textbf{Step 1: Error Decomposition.}
By \Cref{lemma:pathwise-low-noise}, for every $n \geq 1$, it holds that
\begin{align}
\label{eq:thm_low_noise_eq1}
\sum_{k=1}^n \Delta(x_k) &\leq \frac{2}{(2 - L \cdot \eta) \cdot \eta} \cdot \|x_1-x_\star\|^2 + \frac{4}{2-L \cdot \eta} \cdot M_n + G_n \nonumber \\
&\qquad + \frac{2\cdot \eta \cdot (2 + L\cdot \eta)}{(2 - L \cdot \eta)^2} \cdot \sum_{k=1}^n \|\nabla f(x_\star,Z_k)\|^2.
\end{align}
%

\medskip
\noindent\textbf{Step 2: Localization.} Define
\begin{align*}
R_\eta \coloneqq 2 \cdot \|x_1-x_\star\| + 6 \cdot \eta \cdot \sigma_\star \sqrt{ \frac{2 + L \cdot \eta}{2 - L \cdot \eta} \cdot n \cdot \log \Big( \frac{3}{\beta} \Big)}.
\end{align*}
We now show that $\bbP(\tau_{R_\eta}>n) \geq 1 - (2\beta)/3$. 

\noindent Indeed, for every $k\ge1$,
\begin{align*}
\|x_{k+1}-x_\star\|^2 = \|x_k-x_\star\|^2 -2\eta\left\langle \nabla f(x_k,Z_k),x_k-x_\star\right\rangle + \eta^2\|\nabla f(x_k,Z_k)\|^2.
\end{align*}
Using eq.~\eqref{eq:lemma_decomposition_eq2} with $\lambda = ( 2 - L \cdot \eta)/(2 \cdot L \cdot \eta)$ along with the convexity of $f(\cdot, Z_k)$ (\Cref{ass:convexity}), yields
\begin{align*}
\|x_{k+1}-x_\star\|^2 &\leq \|x_k-x_\star\|^2 - \frac{\eta(2-L\eta)}{2} \Big( \underbrace{\Delta(x_k,Z_k) - \langle \nabla f(x_\star,Z_k),x_k-x_\star\rangle}_{H_k} \Big) \\
&\qquad + 2\eta \left\langle \nabla f(x_\star,Z_k),x_\star-x_k \right\rangle + \eta^2\frac{2+L\eta}{2-L\eta}\|\nabla f(x_\star,Z_k)\|^2.
\end{align*}
Since $\eta\in(0,2/L)$ and $f(\cdot, Z_k)$ is convex (\Cref{ass:convexity}), we have that $H_k$ is nonnegative. Hence, summing from $k=1$ to $(m \wedge \tau_{R_\eta})-1$ and dropping the negative term yields, for every $m\in[n]$,
\begin{align*}
\|x_{m\wedge\tau_{R_\eta}}-x_\star\|^2 &\leq \|x_1-x_\star\|^2 + 2\eta M_{(m \wedge \tau_{R_{\eta}})-1}^{(R_{\eta})} + \eta^2\frac{2+L\eta}{2-L\eta} \sum_{k=1}^n \|\nabla f(x_\star,Z_k)\|^2.
\end{align*}
Therefore,
\begin{align*}
\max_{m\in[n]}\|x_{m\wedge\tau_{R_\eta}}-x_\star\|^2 &\leq \|x_1-x_\star\|^2 + 2\eta \max_{m\in[n]} M_{m}^{(R_{\eta})} + \eta^2\frac{2+L\eta}{2-L\eta} \sum_{k=1}^n \|\nabla f(x_\star,Z_k)\|^2.
\end{align*}

\noindent Now define the event
\begin{align*}
E_\eta &\coloneqq \left\{ \max_{m\in[n]} M_{m}^{(R_{\eta})} \le
2 \sigma_\star R_\eta\sqrt{2\log(2)\,n\log(3/\beta)} \right\} \\
&\quad \cap \left\{ \sum_{k=1}^n \|\nabla f(x_\star,Z_k)\|^2 \le \sigma_\star^2\bigl(n\log(2)+\log(3/\beta)\bigr) \right\},    
\end{align*}
and note that, by \Cref{lem:localized_Mn_bound,lem:quadratic_noise_bound}, $\bbP(E_\eta)\ge 1-\frac{2\beta}{3}$. On the event $E_\eta$, we obtain
\begin{align*}
\max_{m\in[n]}\|x_{m\wedge\tau_{R_\eta}}-x_\star\|^2 &\leq \|x_1-x_\star\|^2 + 4\eta \sigma_\star R_\eta\sqrt{2\log(2)\,n\log(3/\beta)} \\
&\qquad + \eta^2\frac{2+L\eta}{2-L\eta} \sigma_\star^2\bigl(n\log(2)+\log(3/\beta)\bigr) \\
&\leq \|x_1-x_\star\|^2 + 4 \cdot \eta \sigma_\star R_\eta\sqrt{2\log(2)\,n\log(3/\beta)} \\
&\qquad + 2 \cdot \eta^2\frac{2+L\eta}{2-L\eta} \sigma_\star^2 n \log(3/\beta).
\end{align*}
Now set
\[
r \coloneqq \|x_1-x_\star\|, \quad t \coloneqq \eta \sigma_\star \sqrt{ \frac{2+L\eta}{2-L\eta} n \log(3/\beta)}.
\]
Then the previous display becomes
\[
\sup_{m\in[n]}\|x_{m\wedge\tau_{R_\eta}}-x_\star\|^2 \le r^2 + 4\sqrt{2\log(2)}\,R_\eta t + 2t^2,
\]
and, by definition, $R_\eta=2r+6t$.
Hence
\[
R_\eta^2-r^2-4\sqrt{2\log(2)}\,R_\eta t-2t^2 = 3r^2+\bigl(24-8\sqrt{2\log(2)}\bigr)rt+\bigl(34-24\sqrt{2\log(2)}\bigr)t^2 \geq 0,
\]
which implies
\[
\sup_{m\in[n]}\|x_{m\wedge\tau_{R_\eta}}-x_\star\|^2 \leq r^2 + 4\sqrt{2\log(2)}\,R_\eta t + 2t^2 \le R_\eta^2 \Rightarrow \tau_{R_{\eta}} > n.
\]

\medskip
\noindent\textbf{Step 3: Control of the stochastic terms.} Define the event
\[
\mathcal G_\eta \coloneqq \left\{ G_n^{(R_{\eta})} \le
\frac12 \sum_{k=1}^n \Delta(x_k)\mathbf 1_{\{k<\tau_{R_\eta}\}} + \frac{2}{3}L R_\eta^2\log(3/\beta) \right\},
\]
and note that by \Cref{lem:localized_Gn_bound}, $\bbP(\mathcal G_\eta)\ge 1-\frac{\beta}{3}$. Let $ E \coloneqq E_\eta \cap \mathcal G_\eta$, and note that by the previous step and the union bound, $\bbP(E) \ge 1-\beta$. Recall also that, on the event $E$, we have $\tau_{R_\eta}>n$. Therefore,
$\mathbf 1_{\{k<\tau_{R_\eta}\}} = 1$ for every $k\in[n]$, and the truncated variables become equal to the original quantities. In particular,
\begin{align}
\label{eq:thm_low_noise_step3_Mn}
M_n &= M_n^{(R_{\eta})}\le 2 \sigma_\star R_\eta\sqrt{2\log(2)\,n\log(3/\beta)},
\end{align}
and
\begin{align}
\label{eq:thm_low_noise_step3_quad}
\sum_{k=1}^n \|\nabla f(x_\star,Z_k)\|^2 \le \sigma_\star^2\bigl(n\log(2)+\log(3/\beta)\bigr) \leq 2 \sigma_{\star}^2 n \log(3/\beta).
\end{align}
Moreover,
\begin{align}
\label{eq:thm_low_noise_step3_Gn}
G_n\leq \frac12 \sum_{k=1}^n \Delta(x_k) + \frac{2}{3}L R_\eta^2\log(3/\beta).
\end{align}

\medskip
\noindent\textbf{Step 4: Conclusion.}
On the event $E$, by eq.~\eqref{eq:thm_low_noise_eq1}, \eqref{eq:thm_low_noise_step3_Mn}, \eqref{eq:thm_low_noise_step3_quad}, and \eqref{eq:thm_low_noise_step3_Gn}, we obtain
\begin{align*}
\sum_{k=1}^n \Delta(x_k) &\leq \frac{2}{(2 - L \cdot \eta) \cdot \eta} \cdot \|x_1-x_\star\|^2 + \frac{4}{2-L \cdot \eta} \cdot M_n + G_n  \\
&\qquad +
\frac{2\cdot \eta \cdot (2 + L\cdot \eta)}{(2 - L \cdot \eta)^2} \cdot \sum_{k=1}^n \|\nabla f(x_\star,Z_k)\|^2 \\
&\leq \frac{2}{(2 - L \cdot \eta) \cdot \eta} \cdot \|x_1-x_\star\|^2 + \frac{4}{2-L \cdot \eta} \cdot M_n +
\frac12 \sum_{k=1}^n \Delta(x_k) + \frac{2}{3}L R_\eta^2\log(3/\beta) \\
&\qquad + \frac{2\cdot \eta \cdot (2 + L\cdot \eta)}{(2 - L \cdot \eta)^2} \cdot
\sum_{k=1}^n \|\nabla f(x_\star,Z_k)\|^2.
\end{align*}
Rearranging, and using eq.~\eqref{eq:thm_low_noise_step3_Mn} and \eqref{eq:thm_low_noise_step3_quad}, yields
\begin{align*}
\sum_{k=1}^n \Delta(x_k) &\leq \frac{4}{(2 - L \cdot \eta) \cdot \eta} \cdot \|x_1-x_\star\|^2 + \frac{8}{2-L \cdot \eta} \cdot M_n +
\frac{4}{3}L R_\eta^2\log(3/\beta) \\
&\qquad + \frac{4\cdot \eta \cdot (2 + L\cdot \eta)}{(2 - L \cdot \eta)^2} \cdot \sum_{k=1}^n \|\nabla f(x_\star,Z_k)\|^2 \\
&\leq \frac{4}{(2 - L \cdot \eta) \cdot \eta} \cdot \|x_1-x_\star\|^2 + \frac{16}{2-L \cdot \eta} \cdot \sigma_\star R_\eta\sqrt{2\log(2)\,n\log(3/\beta)} \\
&\qquad + \frac{4}{3}L R_\eta^2\log(3/\beta) + \frac{8 \cdot \eta \cdot (2 + L\cdot \eta)}{(2 - L \cdot \eta)^2} \cdot
\sigma_\star^2 n\log(3/\beta).
\end{align*}
The claim follows by dividing both sides by $n$, and the convexity of $F$ (\Cref{ass:convexity}).

Now choose
\[
\eta
=
\min\left\{
\frac{1}{4L},
\frac{\|x_1-x_\star\|}{32\,\sigma_\star\sqrt{n\log(3/\beta)}}
\right\},
\]
with the convention that the second quantity is $+\infty$ when $\sigma_\star=0$. Then $\eta\le \frac{1}{4L}$, and therefore
\[
\frac{1}{2-L\eta}\le \frac{4}{7},
\qquad
\frac{2+L\eta}{2-L\eta}\le \frac{9}{7},
\qquad
\frac{8(2+L\eta)}{(2-L\eta)^2}\le \frac{288}{49}<6.
\]
Moreover,
\begin{align*}
R_\eta &= 2\|x_1-x_\star\| + 6\eta \sigma_\star \sqrt{ \frac{2+L\eta}{2-L\eta} n\log(3/\beta) } \leq 2\|x_1-x_\star\| + \frac{18}{32}\sqrt{\frac{1}{7}}\|x_1-x_\star\| \\
&< \frac{9}{4}\|x_1-x_\star\|.
\end{align*}
Hence, 
\begin{align*}
\Delta(\bar x_n)
&\leq
\frac{16}{7}\cdot \frac{\|x_1-x_\star\|^2}{\eta n}
+
\frac{27}{4}L\|x_1-x_\star\|^2\frac{\log(3/\beta)}{n} \\
&\qquad
+
\frac{144}{7}\sqrt{2\log(2)}\,
\sigma_\star \|x_1-x_\star\|
\sqrt{\frac{\log(3/\beta)}{n}}
+
6\eta \sigma_\star^2\log(3/\beta).
\end{align*}
Since
\[
\frac{1}{\eta}
\le
4L
+
\frac{32\,\sigma_\star\sqrt{n\log(3/\beta)}}{\|x_1-x_\star\|},
\]
it follows that
\[
\frac{16}{7}\cdot \frac{\|x_1-x_\star\|^2}{\eta n}
\le
\frac{64}{7}L\|x_1-x_\star\|^2\frac{1}{n}
+
\frac{512}{7}\sigma_\star \|x_1-x_\star\|
\sqrt{\frac{\log(3/\beta)}{n}}.
\]
Finally, if $\sigma_\star>0$, then
\[
6\eta \sigma_\star^2\log(3/\beta)
\le
\frac{3}{16}\sigma_\star \|x_1-x_\star\|
\sqrt{\frac{\log(3/\beta)}{n}},
\]
whereas if $\sigma_\star=0$, the left-hand side is equal to zero. Therefore
\begin{align*}
\Delta(\bar x_n)
&\le
\frac{64}{7}L\|x_1-x_\star\|^2\frac{1}{n}
+
\frac{27}{4}L\|x_1-x_\star\|^2\frac{\log(3/\beta)}{n} \\
&\qquad
+
\left(
\frac{512}{7}
+
\frac{144}{7}\sqrt{2\log(2)}
+
\frac{3}{16}
\right)
\sigma_\star \|x_1-x_\star\|
\sqrt{\frac{\log(3/\beta)}{n}}.
\end{align*}
Since
\[
\frac{512}{7} + \frac{144}{7}\sqrt{2\log(2)} + \frac{3}{16} < 98,
\]
we conclude that
\begin{align*}
\Delta(\bar x_n) &\leq \frac{64}{7}L\|x_1-x_\star\|^2\frac{1}{n} +
\frac{27}{4}L\|x_1-x_\star\|^2\frac{\log(3/\beta)}{n} \\
&\qquad + 98 \sigma_\star \|x_1-x_\star\|
\sqrt{\frac{\log(3/\beta)}{n}},
\end{align*}
which proves item~(ii).

\subsection{Nonconvex Optimization}\label{appendix:extensions_nonconvex_optimization}
The result of this section is the following.
\begin{theorem}
\label{thm:pl_qg_extension_full}
Suppose Assumption~\ref{ass:feasibility} holds. Then the iterates generated by eq.~\eqref{eq:sgd_cs} satisfy the following.
\begin{enumerate}[label={\rm(\roman*)}]
\item\label{item:pl_expectation_bound_main} Suppose Assumptions~\ref{ass:interpolation} and \Cref{ass:nonconvexity}~\cref{ass:nonconvexity_PL} also hold. Then, for every $n\geq 1$, it holds that
\begin{align}
\label{eq:pl_expectation_bound_main}
\bbE[\Delta(x_{n+1})] \leq \Delta(x_1) \cdot (1 - (2 \cdot \mu \cdot \eta - L^2 \cdot \eta^2))^n.
\end{align}
\item\label{item:pl_qg_hp_bound_main} Suppose \Cref{ass:nonconvexity} holds. Then, for every $n\geq 1$ and every $\beta\in(0,1)$, it holds that
\begin{align}
\label{eq:pl_qg_hp_bound_main}
\bbP\bigg( \Delta(x_{n+1}) \leq \frac{1}{\sqrt{\beta}} \cdot \Delta(x_1) \cdot ( 1 - \varphi_{\mu, L}(\eta) )^{n/2} \bigg) \geq 1-\beta,
\end{align}
where
\begin{align*}
\varphi_{\mu, L}(\eta) &\coloneqq 4 \cdot \mu \cdot \eta - (4 \cdot \mu^2 + 10 \cdot L^2 ) \cdot \eta^2 + 4 \cdot \mu \cdot L^2 \cdot \eta^3 - L^4 \cdot \eta^4 - 2 \cdot L^5 \cdot \eta^4/\mu.
\end{align*}
\end{enumerate}
\end{theorem}
\begin{remark}
We make the following comments.
\begin{enumerate}
\item We do not claim novelty on \cref{item:pl_expectation_bound_main} of \Cref{thm:pl_qg_extension_full} as it appeared already in prior works (although for the special case of $F$ being a sum and for $\eta_k$ being SPS) \cite{Loizou2021}. Indeed, we report it here for the sake of comparison with the high-probability rate, and we report the proof for completeness. Besides, it is worth noting that, eq.~\eqref{eq:pl_expectation_bound_main} gives a convergence rate for any $\eta \in (0, 2 (\mu/L^2))$. The rate is optimized at $\eta = \mu/L^2$, where it is of the order of $\exp(- (\mu/L)^2 \cdot n)$.

\item For the high-probability result, we note that 
\[
1-\varphi_{\mu, L}(\eta) = (1-2\mu\eta+L^2\eta^2)^2 + 8 L^2 \eta^2 + 2 L^5 \eta^4/\mu \geq 0
\]
so that the rate in \eqref{eq:pl_qg_hp_bound_main} is well defined for every $\eta$. On the other hand, the bound yields convergence only for $\eta$ s.t. $1-\varphi_{\mu, L}(\eta) < 1$. Algebraic calculations show that this is indeed the case of the choice $\eta = (1/5) (\mu/L^2)$ made in \Cref{thm:pl_qg_extension} in the main body.
\end{enumerate}
\end{remark}

\subsubsection{Proof of Theorem \ref{thm:pl_qg_extension_full}}
The following auxiliary lemma (Proposition 1 from \cite{Apidopoulos2022}) 
shows that functions that satisfy PL also enjoy the following quadratic growth condition.
%
\begin{lemma}\label{lemma:PL_to_QG}
Suppose Assumptions~\ref{ass:feasibility} and \ref{ass:nonconvexity}~\cref{ass:nonconvexity_PL}. Let $\text{dist}(x,X_{\star}):= \inf_{x_\star\in X_{\star}} \|x-x_\star\|$. Then $F$ satisfies 
\begin{align}\label{eq:QG}
\Delta(x) \geq \frac{\mu}{2}\cdot \text{dist}(x,X_{\star})^2.
\end{align}
\end{lemma}

\noindent We now provide a proof of the two claims.

\noindent  \rm{(i)}: By \Cref{ass:nonconvexity} $F$ is $L$-smooth. Thus, for every $k\geq 1$,
\begin{align*}
F(x_{k+1}) &\leq F(x_k) + \langle \nabla F(x_k),x_{k+1}-x_k\rangle + \frac{L}{2} \|x_{k+1}-x_k\|^2 \\
&\overset{(a)}{=} F(x_k) - \eta \langle \nabla F(x_k),\nabla f(x_k,Z_k)\rangle + \frac{L \cdot \eta^2}{2} \|\nabla f(x_k,Z_k)\|^2  \\
&\overset{(b)}{\leq} F(x_k) - \eta \langle \nabla F(x_k),\nabla f(x_k,Z_k)\rangle + L^2 \eta^2 \Delta(x_k, Z_k) \\
&= F(x_k) - \eta \|\nabla F(x_k)\|^2 + \eta \langle \nabla F(x_k),\nabla F(x_k) - \nabla f(x_k,Z_k)\rangle + L^2 \eta^2 \Delta(x_k, Z_k) \\
&\overset{(c)}{\leq} F(x_k) - 2 \eta \mu \Delta(x_k) + \eta \langle \nabla F(x_k),\nabla F(x_k) - \nabla f(x_k,Z_k)\rangle + L^2 \eta^2 \Delta(x_k, Z_k),
\end{align*}
where: (a) follows from the update rule \eqref{eq:sgd_ls}, (b) follows from the smoothness of $f(\cdot, Z_k)$ (\Cref{ass:nonconvexity}) and the fact that $x_{\star}$ is one of its minimizers (\Cref{ass:interpolation}), and (c) follows from $F$ being $\mu$-PL (\Cref{ass:nonconvexity} \cref{ass:nonconvexity_PL}). The above implies
\begin{align}
\label{eq:pl_pw1}
\Delta(x_{k+1}) &\leq (1 - 2 \eta \mu + L^2 \eta^2) \Delta(x_k) + \eta \langle \nabla F(x_k),\nabla F(x_k) - \nabla f(x_k,Z_k)\rangle \nonumber\\ 
&\qquad + L^2 \eta^2 (\Delta(x_k, Z_k) - \Delta(x_k)).
\end{align}
Taking conditional expectation on both sides of the above, and noting that $\nabla F(x_k)$ is $\sigma(Z_1, \dots, Z_{k-1})$-measurable, yields
\begin{align*}
\bbE_k[\Delta(x_{k+1})] &\leq (1 - 2 \mu \eta + L^2 \eta^2) \Delta(x_k).
\end{align*}
Finally, taking total expectation and unrolling the above recursion yields
\begin{align*}
\bbE[\Delta(x_{n+1})] \leq (1 - 2 \mu \eta + L^2 \eta^2)^n \Delta(x_1),    
\end{align*}
which proves eq.~\eqref{eq:pl_expectation_bound_main}.

\medskip

\noindent\rm{(ii)}: From eq.~\eqref{eq:pl_pw1} let us define, for every $k \geq 1$, the following
\begin{align*}
M_k &\coloneqq \eta \cdot \langle \nabla F(x_k), \nabla F(x_k)-\nabla f(x_k,Z_k)\rangle,\\
N_k &\coloneqq L^2 \cdot \eta^2 \cdot \big(\Delta(x_k,Z_k)-\Delta(x_k)\big),
\end{align*}
and note that $\bbE_k[M_k] = \bbE_k[N_k] = 0$. Moreover,
\begin{align}
\label{eq:pl_Mk_var}
\bbE_k[M_k^2] &= \bbE_k[\eta^2 \cdot (\langle \nabla F(x_k), \nabla F(x_k)-\nabla f(x_k,Z_k)\rangle)^2] \nonumber \\
&\overset{(a)}{\leq} \eta^2 \cdot \bbE_k[\|\nabla F(x_k) \|^2 \cdot  \| \nabla F(x_k)-\nabla f(x_k,Z_k)\|^2] \nonumber \\
&\overset{(b)}{=} \eta^2 \cdot \|\nabla F(x_k) \|^2 \cdot  \bbE_k[ \| \nabla F(x_k)-\nabla f(x_k,Z_k)\|^2] \nonumber \\
&\leq \eta^2 \cdot \|\nabla F(x_k) \|^2 \cdot  \bbE_k[ \| \nabla f(x_k,Z_k)\|^2] \nonumber \\
&\overset{(c)}{\leq} 4 \cdot L^2 \cdot \eta^2 \cdot \Delta(x_k) \cdot  \bbE_k[ \Delta(x_k, Z_k)] \nonumber \\
&= 4 \cdot L^2 \cdot \eta^2 \cdot \Delta(x_k)^2,
\end{align}
where: (a) follows from Cauchy--Schwarz inequality, (b) is due to $\|\nabla F(x_k)\|^2$ being $\sigma(Z_1,\dots,Z_{k-1})$-measurable, and (c) follows from $F$ and $f(\cdot, Z_k)$ being $L$-smooth (\Cref{ass:nonconvexity}). Similarly it holds that
\begin{align}
\label{eq:pl_Nk_var}
\bbE_k[N_k^2] &= L^4 \cdot \eta^4 \cdot \bbE_k\big[(\Delta(x_k,Z_k)-\Delta(x_k))^2\big] \nonumber \\
&\leq L^4 \cdot \eta^4 \cdot \bbE_k\big[\Delta(x_k,Z_k)^2\big] \nonumber \\
&\leq L^4 \cdot \eta^4 \cdot \bbE_k\big[\Delta(x_k,Z_k) \cdot  \Delta(x_k, Z_k)\big] \nonumber \\
&\overset{(a)}{\leq} \frac{L^5 \cdot \eta^4}{2} \cdot \bbE_k\big[\text{dist}(x_k, X_{\star})^2 \cdot  \Delta(x_k, Z_k)\big] \nonumber \\ 
&\overset{(b)}{\leq} \frac{L^5 \cdot \eta^4}{\mu} \cdot \bbE_k\big[\Delta(x_k) \cdot  \Delta(x_k, Z_k)\big] \nonumber \\
&= \frac{L^5 \cdot \eta^4}{\mu} \cdot \Delta(x_k)^2,
\end{align}
where (a) follows from $f(\cdot, Z_k)$ being $L$-smooth and \Cref{ass:nonconvexity}, and (b) follows from \Cref{lemma:PL_to_QG}.

\noindent Now, notice that, using the definition of $M_k$ and $N_k$, eq.~\eqref{eq:pl_pw1} can be rewritten as 
\begin{align*}
\Delta(x_{k+1}) \leq (1 - 2 \cdot \mu \cdot \eta + L^2\cdot \eta^2 ) \cdot \Delta(x_k)+M_k+N_k,
\end{align*}
and since the right-hand side is non-negative, squaring both sides and taking conditional expectation yields
\begin{align*}
\bbE_k[\Delta(x_{k+1})^2] &\leq \bbE_k\Big[((1 - 2 \cdot \mu \cdot \eta + L^2\cdot \eta^2 ) \cdot \Delta(x_k)+M_k+N_k)^2\Big]\\
&= (1 - 2 \cdot \mu \cdot \eta + L^2\cdot \eta^2 )^2 \cdot \Delta(x_k)^2 \\
&\quad + 2 \cdot (1 - 2 \cdot \mu \cdot \eta + L^2\cdot \eta^2 ) \cdot \Delta(x_k) \cdot \bbE_k[M_k+N_k] + \bbE_k[(M_k+N_k)^2] \\
&\overset{(a)}{\leq} (1 - 2 \cdot \mu \cdot \eta + L^2\cdot \eta^2 )^2 \cdot \Delta(x_k)^2 + 2 \cdot (\bbE_k[M_k^2] + \bbE_k[N_k^2]) \\
&\overset{(b)}{\leq} (1 - 2 \cdot \mu \cdot \eta + L^2\cdot \eta^2 )^2 \cdot \Delta(x_k)^2 + 2 \cdot (4 \cdot L^2 \cdot \eta^2 + L^5 \cdot \eta^4/\mu) \cdot \Delta(x_k)^2 \\
&= \left((1 - 2 \cdot \mu \cdot \eta + L^2\cdot \eta^2 )^2 + 8 \cdot L^2 \cdot \eta^2 + 2 \cdot L^5 \cdot \eta^4/\mu \right) \cdot \Delta(x_k)^2 \\
&= (1 - \varphi_{\mu, L}(\eta)) \cdot \Delta(x_k)^2,
\end{align*}
where
\begin{align*}
\varphi_{\mu, L}(\eta) \coloneqq 4 \mu \eta - (4\mu^2 + 10L^2 ) \eta^2 + 4 \mu L^2 \eta^3 - L^4 \eta^4 - 2 L^5 \eta^4/\mu,
\end{align*}
and where $(a)$ follows from the fact that $M_k$ and $N_k$ have (conditionally) 0-mean, and $(b)$ follows from eq.~\eqref{eq:pl_Mk_var} and eq.~\eqref{eq:pl_Nk_var}. 
Taking total expectation and unrolling yields
\begin{align}
\label{eq:pl_qg_second_moment_bound_q2}
\bbE[\Delta(x_{n+1})^2] \leq ( 1 - \varphi_{\mu, L}(\eta) )^n \cdot \Delta(x_1)^2.
\end{align}
Finally, define
\begin{align*}
t \coloneqq \frac{1}{\sqrt{\beta}} \cdot ( 1 - \varphi_{\mu, L}(\eta) )^{n/2} \cdot \Delta(x_1).
\end{align*}
By Markov's inequality and eq.~\eqref{eq:pl_qg_second_moment_bound_q2},
\begin{align*}
\bbP(\Delta(x_{n+1})\geq t) = \bbP(\Delta(x_{n+1})^2\geq t^2) \leq \frac{\bbE[\Delta(x_{n+1})^2]}{t^2} \leq \beta,
\end{align*}
or equivalently,
\[
\bbP(\Delta(x_{n+1})\leq t)\geq 1-\beta,
\]
which proves eq.~\eqref{eq:pl_qg_hp_bound_main}.

\end{document}